# Multi-Objective Multi-mode Time-Cost Tradeoff modeling in Construction Projects Considering Productivity Improvement

Ali Mohammadjafari [a], Seyed Farid Ghannadpour [b], Morteza Bagherpour [c], Fatemeh Zandieh [d]

[a] Department of Computer Science, University of Louisiana at Lafayette, USA.

Email: ali.mohammadjafari1@louisiana.edu

[b] Department of Industrial Engineering, Iran University of Science and Technology, Iran, 16846-13114.

Email: ghannadpour@iust.ac.ir

Tel: (+9821)-73225015

[c] Department of Industrial Engineering, Iran University of Science and Technology, Iran, 16846-13114.

Email: bagherpour@iust.ac.ir

[d] Department of Industrial Engineering, Iran University of Science and Technology, Iran, 16846-13114.

Email: fatemeh_zandiyeh@ind.iust.ac.ir

Corresponding author (✉)

Permanent Address: Department of Industrial Engineering, Iran University of Science and Technology, Tehran, Iran, 16846-13114.

Tel.: +98 21 73225015;

Email: ghannadpour@iust.ac.ir

# Multi-Objective Multi-mode Time-Cost Tradeoff modeling in Construction Projects Considering Productivity Improvement

**Abstract**

In today's construction industry, poor performance often arises due to various factors related to time, finances, and quality. These factors frequently lead to project delays and resource losses, particularly in terms of financial resources. This research addresses the Multimode Resource-Constrained Project Scheduling Problem (MRCPSP), a real-world challenge that takes into account the time value of money and project payment planning. In this context, project activities exhibit discrete cost profiles under different execution conditions and can be carried out in multiple ways. This paper aims to achieve two primary objectives: minimizing the net present value of project costs and project completion times while simultaneously improving the project's productivity index. To accomplish this, a mathematical programming model based on certain assumptions is proposed. Several test cases are designed, and they are rigorously evaluated using the methodology outlined in this paper to validate the modeling approach. Recognizing the NP-hard nature of this problem, a multi-objective genetic algorithm capable of solving large-scale instances is developed. Finally, the effectiveness of the proposed solution is assessed by comparing it to the performance of the NSGA-II algorithm using well-established efficiency metrics. Results demonstrate the superior performance of the algorithm introduced in this study.

## 1. Introduction

Multi-mode project scheduling problems, while striving to enhance productivity, are prevalent challenges in real-world industries. Numerous endeavors have been dedicated to addressing multi-mode project scheduling issues with the aim of improving productivity. However, the integrated consideration of multi-mode scheduling alongside multi-skill allocation problems has been relatively rare.

The construction industry stands as one of humanity's most crucial sectors. Its significance has been consistently evident in societies throughout history and continues to grow in tandem with population expansion. Prior to the COVID-19 pandemic, the construction industry had already achieved an impressive market value of approximately $12 trillion, with expectations of annual growth around 3%. According to global statistics, it is forecasted that construction costs will soar to $14,000 billion by 2025.

Nonetheless, underperformance remains a pervasive concern within the construction sector, stemming from temporal, cost, and quality-related factors. These issues often lead to project delays and substantial resource losses, particularly in terms of finances. Official statistics illustrate this challenge vividly: in the past three years, only 25% of projects managed to adhere to delivery schedules with a delay of no more than 10% of the initially contracted time. Similarly, a mere 31% of projects were completed within budgets, exceeding predictions by no more than 10%. Large-scale projects tend to extend their timelines by up to 20% beyond the initially estimated duration, accompanied by costs that can surpass the original budget by as much as 80%. These reports underscore the pressing need for research that can introduce a methodology capable of effectively balancing project delivery times and overall costs, all while considering the broader context of total project productivity.

The Multimode Resource-Constrained Project Scheduling Problem (MRCPSP) represents a fundamental challenge in project scheduling, extending beyond the classical Resource-Constrained Project Scheduling Problem (RCPSP). MRCPSP models real-world scenarios by accounting for multiple execution modes, each associated with distinct costs, resources, and durations for a given activity. Unlike traditional scheduling, where activity durations remain constant, MRCPSP acknowledges that durations can vary depending on the chosen execution mode. These variations may involve incurring additional costs for faster completion or reductions in expenses for longer durations. As a result, this study introduces discrete and varying cost slopes for each activity. Additionally, the Project Payment Plan (PPP) plays a pivotal role in project cash flow management, offering four distinct approaches for its implementation (Ulusoy et al., 2001):

✓ Lump-Sum Payment (LSP): The total cost is paid to the contractor when the project terminates successfully.
✓ Payments at Event Occurrence (PEO): Payments are made when events occur.
✓ Payments at Activities' Completion time (PAC): It is made after the completion of every activity.
✓ Progress Payments (PP): Payments are made in regular time intervals, and the last payment is made when the project is completed.

In this study, the Payment Event Occurrence (PEO) model is employed. Under this model, a predefined set of project activities are designated as significant payment milestones. When these specific activities are completed, a predetermined portion of the contract amount is disbursed to the contractor, in accordance with the corresponding formulas provided in equations (1) through (3). These formulas serve to ascertain the payment events:

$$i_j = \left\{ i: y_{iT} = 1, T = min\left[ T: \sum_{i=1}^{n} \left( v_i \sum_{t=0}^{T} y_{it} \right) \geq j\left(\frac{U}{J}\right) \right] \right\} \qquad j = 1,..J-1 \qquad (1)$$

$$i_j = \left\{ i: y_{iT} = 1, T = min\left[ T: \sum_{i=1}^{n} \left\{ \left[ \sum_{m=1}^{M} (C_{im} X_{im}) \right] \sum_{t=0}^{T} y_{it} \right\} \geq j\left(\frac{B}{J}\right) \right] \right\} \qquad j = 1,..J-1 \qquad (2)$$

$$i_j = \left\{ i: y_{iT} = 1, T = min\left[ T: T \geq j\left(\frac{D}{J}\right) \right] \right\} \qquad j = 1,..J-1 \qquad (3)$$

Formula (1) requires that the j-th payment corresponds to the earliest activity that brings the contractor's cumulative earned value to a point where it equals or exceeds j (U/J). Formula (2) mandates that the j-th payment is associated with the earliest activity that causes the contractor's cumulative expenses to reach or surpass j (B/J), with B representing the project's benchmark cost for determining payment milestones. Formula (3) guarantees that the j-th payment is allocated to the earliest activity that completes no sooner than j (D/J). For the purposes of this paper, Formula (3) is employed to determine payment milestones.

In the contemporary competitive landscape, productivity improvement is paramount for the survival and success of organizations and projects. It serves as a critical factor in project success, encompassing aspects like product quality and cost efficiency. Productivity is highly valued across various projects, as it enables contractors to maximize profits while minimizing time and expenses through the use of productivity enhancement methods. Neglecting productivity can lead to project disruptions and customer dissatisfaction. Thus, enhancing project quality is a key approach to boosting overall project productivity.

This study delves into a three-objective Multimode Resource-Constrained Project Scheduling Problem (MRCPSP), incorporating considerations for the discounted time value of money and Payment Event Occurrence (PEO)-based payments. Notably, this problem accounts for the allocation of both renewable and non-renewable resources to project activities. The primary objective is to calculate and minimize the net present value of project costs by factoring in the time value of money. The secondary objective aims to reduce project completion time, specifically based on the last event occurrence. The third objective focuses on enhancing overall project productivity, measured as the ratio of project quality output to total (direct and indirect) project costs. In this research, each activity is assigned a discrete cost slope that varies depending on its execution mode. The ε-constraint method is employed to balance these three objectives, serving as a validation approach for the proposed model. This model is subsequently implemented and solved using the BARON solver within the GAMS software. Given the NP-hard nature of the problem, two efficient meta-heuristic Pareto-based algorithms, the Non-Dominated Sorting Genetic Algorithm (NSGA-II) and multi-objective Genetic Algorithm (MOGA), are developed for its resolution. In light of the problem's significance, this paper's contributions are outlined below:

1. Developing a three-objective mathematical model to balance the completion time and the net present value of the project's total costs, aiming to enhance productivity.
2. Incorporating the time value of money and cash flow into the project, while also accommodating project payments through the PEO method.
3. Assessing the quality index of project activities using the DANP method and calculating productivity based on this index.
4. Adapting the project to a multimode condition, where each activity exhibits distinct and discrete cost slopes in every mode.

The organization of the subsequent sections of the paper is as follows: Section 2 reviews the research literature, Section 3 states and describes the proposed problem and mathematical model, Section 4 tackles the solving method of the problem, and Section 5 explains the results of solving small- and large size problems by the two algorithms and compares their outputs. Finally, Section 6 encompasses the research conclusion and vision.

## 2. Literature review

### 2.1. Resource-Constrained Project Scheduling Problem (RCPSP)

This section of the literature review discusses the project scheduling problem, particularly focusing on the role of constrained resources and their impact on project costs and cash flows. It highlights several key research contributions in this domain.

Liu and Wang (2008) introduced a project planning model integrating resource and cash flow considerations, with a focus on optimizing contractor profits. Chen and Weng (2009) proposed a two-phase genetic algorithm (GA) that addressed time-cost tradeoffs and resource allocation in project scheduling.

Chen and Shahandashti (2009) and Niño and Peña (2019) improved upon this problem with a simulated annealing-genetic algorithm (SA-GA) hybrid and a honeybee algorithm, respectively. Rahman et al. (2020) presented a GA-based Memetic Algorithm (MA) that demonstrated its efficacy through numerical results and comparisons with advanced algorithms. To evaluate scheduling policies under uncertainty, Rostami et al. (2018) employed the Markov Chain Monte Carlo method for the Stochastic Resource-Constrained Project Scheduling Problem (SRCPSP). This method, although time-consuming, enhanced solution quality significantly.

Roghanian et al. (2018) introduced a modified critical chain approach using fuzzy logic to minimize project implementation time under uncertain conditions. Birjandi and Mousavi (2019) developed a method based on fuzzy mixed-integer non-linear programming (MINLP) for selecting scheduling paths in uncertain environments. Su et al. (2020) introduced a theory for efficient time conservation period allocation to prevent management delays. Cui et al. (2021) addressed resource-constrained multi-project scheduling with an integrated multi-mode and multi-skill model for high-end equipment development. Yuan et al. (2021) tackled prefabricated building (PB) construction scheduling with a hybrid cooperative co-evolution algorithm to handle execution time uncertainty and multiple objectives. Liu et al. (2022) considered stochastic duration and resource demand in a chance constrained programming model to minimize project duration and resource costs. Sayyadi et al. (2022) developed a bi-objective model that combines Resource Leveling Problem (RLP) and Multi-Project Scheduling Problem (MPSP) to minimize project durations and resource usage, integrating a community detection approach. Milička et al. (2022) studied a multi-agent project employee problem under constrained resources, proposing a bi-level optimization model considering interactions between team leaders and project managers. Additionally, more studies in this area can be found in the works of Pass-Lanneau et al. (2023), Van Eynde et al. (2023), He et al. (2023), Liu et al. (2023), and Phuntsho and Gonsalves (2023).

### 2.2. Productivity and Efficiency

This section explores the topic of productivity and efficiency in project scheduling, with a focus on various models and approaches proposed by different researchers.

Mehmanchi and Shadrokh (2013) and Zabihi et al. (2019) introduced models for addressing the multi-skill project scheduling problem (MSPSP). They incorporated a dynamic exponential learning function, considering employee efficiency assumptions and the learning effect. To linearize this non-linear function, they employed separable programming. Ahmadpour and Ghezavati (2019) presented a fuzzy scheduling model for MSPSP, which involved work calendars for project members and skill factor determination based on efficiency concepts. Lin and Chou (2019) developed a meta-heuristic genetic algorithm aimed at minimizing project time. Stylianou and Andreou (2016) used multi-objective optimization to simultaneously control project cost and time. They considered various productivity-related factors to estimate activity costs and durations, including developer efficiency and extra relational costs. Jeunet and Bou Orm (2020) explored the relationship between quality and human resources, considering factors like human force efficiency and its impact on quality due to overtime work. They aimed to optimize the number of permanent and temporary workers and overtime work to reduce time, cost, and enhance personal quality. Wang et al. (2020) introduced a multi-skill scheduling programming approach that allowed for the interruption of main activities to enhance productivity and efficiency. They considered factors like employee fatigue over time and the impact of shifts between different jobs on employee weariness and productivity. These studies collectively contribute to the understanding of productivity and efficiency in project scheduling, offering various models and techniques to address these critical aspects.

## 2.3. Multi-Objective Multi-Mode RCPSP

Wang and Zheng (2018) proposed a multi-objective Fruit Fly Optimization algorithm for the Multi-Skill Resource-Constrained Project Scheduling Problem (MSRCPSP), aiming to minimize both time and costs concurrently. Balouka and Cohen (2021) introduced a Robust Optimization (RO) approach to solve the Multimode Resource-Constrained Project Scheduling Problem under uncertain activity durations. They employed the Benders decomposition method to jointly consider mode selection and resource allocation, ultimately reducing resource allocation. Maghsoudlou et al. (2017) investigated the original version of the multi-skill Resource-Constrained Project Scheduling Problem, focusing on the duplication risk of activities based on multi-skill human resource allocation. They proposed a two-objective optimization model to minimize activity processing costs and duplication risks simultaneously. Tirkolaee et al. (2019) conducted an MMRCPSP study with a focus on maximizing the net present value and minimizing project completion time, considering renewable resources such as human labor. They developed the NSGA-II and Multi-Objective Simulated Annealing (MOSA) algorithms and found MOSA to be more efficient for small-size problems, while NSGA-II performed better for large-size problems. Fernandes Muritiba et al. (2018) introduced the Path-Rethinking (PR) algorithm for solving the MRCPSP. Eydi and Bakhshi (2019) presented a multi-objective model for solving an RCPSP problem, aiming to minimize project time while maximizing the net present value of project cash flows. Nemati-Lafmejani et al. (2019) proposed an uncertain two-objective optimization model for solving the Multimode Resource-Constrained Project Scheduling and Contractor Selection Problems (MRCPSP-CS). They conducted sensitivity analysis, showing that increasing the number of contractors could offer more flexible decision-making options. Subulan (2020) presented an approach based on interval-stochastic programming for the MRCPSP, considering various project scheduling risks, human resources, and achieving a trade-off between project time and total human resource costs. Fernandes and de Souza (2021) addressed the MRCPSP by considering activity cessation conditions. They aimed to minimize project time and identify the start of each project activity, proposing a local bifurcation mathematical strategy for solving the problem. Józefowska et al. (2001) and Sadeghloo et al. (2023) are other relevant studies in this field.

These studies contribute to the understanding and solutions for the Multi-Objective Multi-Mode Resource-Constrained Project Scheduling Problem, offering various models and algorithms to handle its complexity.

## 2.4. Time-Cost Tradeoff Problem of Projects

This section delves into the Time-Cost Tradeoff (TCT) problem in the context of project management, which involves optimizing resource allocation to strike a balance between competing project aspects. The TCT problem has been a significant focus for construction engineering and management researchers, particularly in minimizing project time and cost while addressing various challenges.

Ghoddousi et al. (2013) considered the Multimode Resource-Constrained Project Scheduling Problem (MRCPSP), Discrete Time-Cost Tradeoff Problem (DTCTP), Resource Allocation, and Resource Leveling Problem (RLP) simultaneously. Ke and Ma (2014) tackled complex environments with multiple uncertain conditions, incorporating fuzzy stochastic theory to describe the project environment. Nabipoor Afruzi et al. (2014) explored the discrete time-cost-quality tradeoff problem (DTCQTP) in multimode resource-constrained projects, addressing real-world resource constraints. Tavana et al. (2014) proposed a new multimode multi-objective model with preconditions relations, aligning it with real-world project scenarios. He et al. (2017) investigated a discrete time-cost tradeoff problem while minimizing the maximal cash flow gap, a significant metric for contractor cash flows during a project. Leyman et al. (2019) studied different solutions' effects on optimizing the net present value in project planning, particularly in the context of the time-cost tradeoff problem. Tareghian and Taheri (2007) developed a simulation-based integer linear programming tool to assess project feasibility and profitability in the presence of risks. Amoozad Mahdiraji et al. (2021) proposed an approach to identify the best implementation situation for each project activity, considering time, cost, quality, and risk criteria under uncertain circumstances. Panwar and Jha (2021) developed a decision-making model that incorporated time, cost, quality, and safety (TCQS) at the planning stage, utilizing a many-objective evolutionary algorithm.

In addition to addressing time and cost, this study conducted a correlation analysis to explore relationships among the TCQS components. Table 1 in the paper presents a research gap analysis based on the literature, highlighting areas where further investigation and development are needed.

**Table 1. Literature overview**

| Reference | Problem's Goals | | | | MADM | TVM | RR | MP | Multimode | DDCS | Solution methods |
|---|---|---|---|---|---|---|---|---|---|---|---|
| | Cost | Quality | Time | Productivity | | | | | | | |
| Tareghian and Taheri (2007) | * | * | * | | | | * | | | | Metaheuristics |
| Liu and Wang (2008) | * | * | | | | * | * | | | | Constraint Programming |
| Ghoddousi et al. (2013) | * | | * | | | | * | | * | | NSGA-II |
| Khalili-Damghani et al. (2015) | * | * | * | | | | * | | * | | Classical epsilon-constraint |
| Tran et al. (2015) | * | * | * | | | | * | | | | MOABCDE-TCQ |
| Stylianou and Andreou (2016) | * | | * | * | | | * | | | | MOCell, NSGA-II , PAES, ,and SPEA2 |
| Rostami et al. (2018) | | | * | | | | * | | | | Two phase local search algorithm |
| da Silva et al. (2017) | * | | * | * | | | * | | | | DEA |
| Tirkolaee et al. (2019) | | * | * | | | | * | * | * | | NSGA-II + MOSA |
| Balouka and Cohen (2021) | | | * | | | | | | * | | Robust Optimization |
| Eydi and Bakhshi (2019) | | * | * | | | * | * | | | | NSGA-II |
| Zabihi et al. (2019) | | | * | * | | | * | | | | TLBOT |
| Luong et al. (2021) | * | * | * | | | | * | | * | | MOPSO, NSGA-II |
| Subulan (2020) | * | | * | | | | * | | * | | Hybrid interval programming |
| Rahman et al. (2020) | | | * | | | | * | | | | GA-based memetic algorithm (MA) |
| Fernandes and de Souza (2021) | | | * | | | | * | | * | | Local Branching |
| Elloumi et al. (2021) | | | * | | | | * | | * | | Hurried Slipping Window Method |
| This paper | * | * | * | * | * | * | * | * | * | * | MOGA |

**MADM:** Multi-Attribute Decision Making, **TVM:** Time Value of Money, **RR:** Renewable Resource, **MP:** Method of Payment, **DDCS:** Different Discrete Cost Slope.

# 3. Multi-mode Productivity-based Time-Cost Tradeoff Model

Within the realm of project management, the time-cost tradeoff problem, particularly in construction projects, has recently gained significant theoretical and practical attention. This problem revolves around the challenge of accelerating a project's schedule by expediting certain activities and allocating additional resources. While this can lead to higher direct costs, it can also reduce indirect expenses like overhead, fixed costs, and daily expenses by shortening the project's duration. The intriguing aspect is that the outcome of schedule compression on project costs can vary, making it possible to simultaneously reduce both time and expenses. In this study, each project activity can assume multiple execution modes, each associated with different timeframes and costs, contingent upon the allocated resources. Time-cost tradeoff optimization is a process aimed at identifying effective methods for executing project activities, with the goal of accelerating the project and determining the most suitable time-cost mode. This process involves selecting the most appropriate resources, such as labor, equipment, methods, and technology, for each project activity. Given the multitude of potential combinations for these resource alternatives, the challenge lies in determining the optimal combination to achieve the best balance between time and cost while adhering to the project's unique time and budget constraints. Figure1 provides a schematic representation of the research problem.

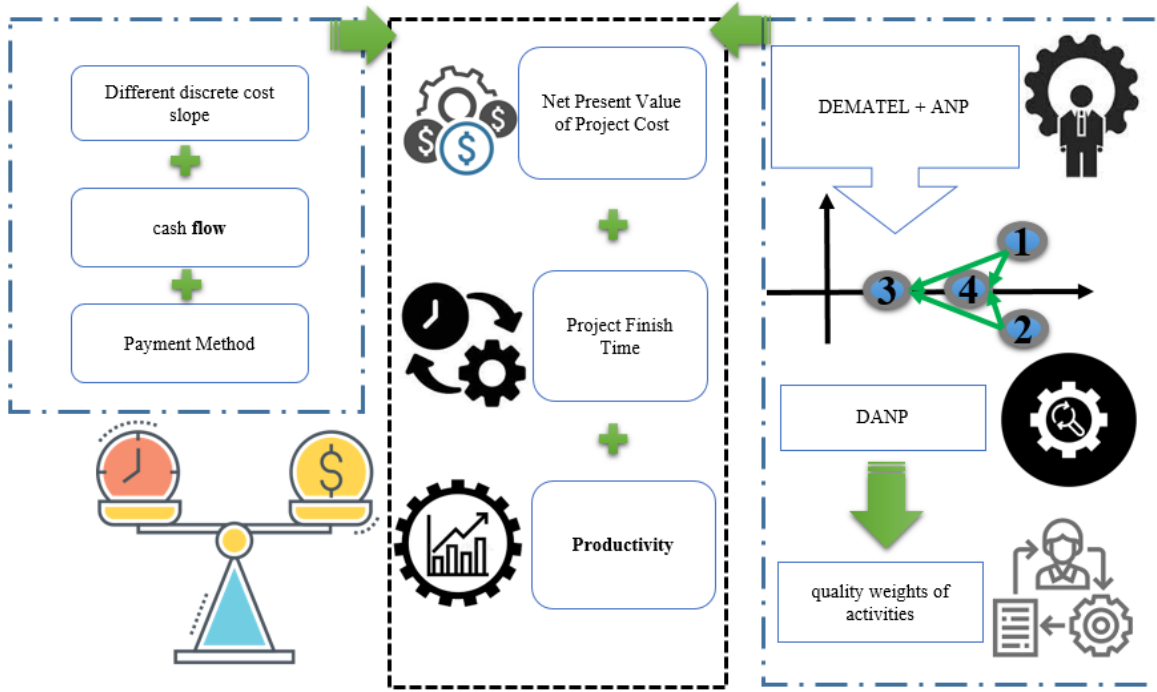

**Figure 1. Schematic of the Research Problem**

## 3.1. Different Discrete Cost Slopes in Multimode Project Scheduling

To address real-world challenges effectively, this study introduces a novel concept: different and discrete cost slopes for each execution mode of a project's activities. Essentially, it considers time and cost as functions of the compression rate, which signifies the reduction of activity duration at the expense of increased direct costs. For instance, as illustrated in Figure, in a given project, Worker A can complete Activity $i$ using Machine X (mode 1) within 10 working days, each consisting of 8 hours, at a cost of $400. On the other hand, Worker B can opt for Machine Y (mode 2) to complete the same activity in 8 working days, also with 8-hour shifts, but at a cost of $500. Importantly, both Worker A and Worker B can expedite the activity by working overtime, which inevitably raises costs. In essence, time shortening denotes the act of accelerating activity completion through various means, such as overtime work or the allocation of additional resources, which may incur extra expenses but decrease project duration. Consequently, the cost of executing Activity $i$ is not solely contingent upon the chosen mode but also influenced by the selected duration.

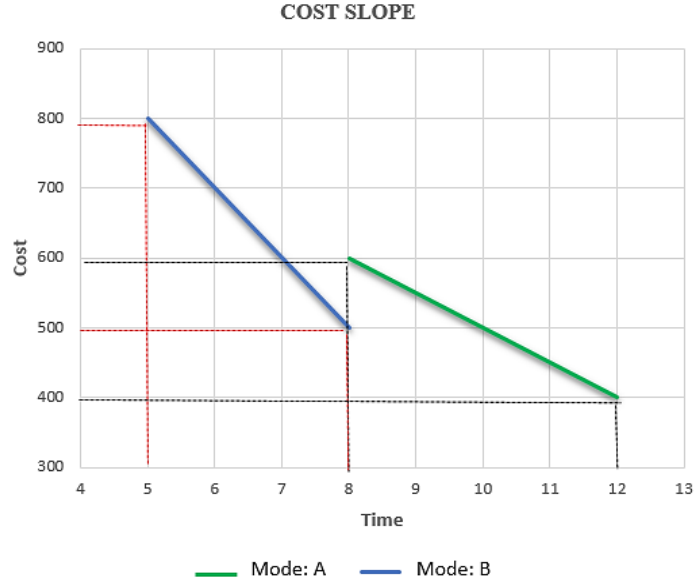

**Figure 2. Different Discrete Cost Slopes**

The central goal of this research is to address the time-cost tradeoff problem, which involves a comprehensive analysis of project costs in relation to variations in the timeframes for activity execution. The aim is to identify the optimal composition for reducing activity durations in a way that minimizes the total project costs. In practical project scenarios, where projects often span extended time periods, factors like the time value of money and interest rates significantly impact overall project costs. Effectively managing and mitigating this impact over the long term is a fundamental concern for contractors. Therefore, this study endeavors to reduce the net present value of project costs by implementing a Payment at Event Occurrence model within the Project Payment Plan. The tradeoff model developed for this research, drawing inspiration from the models proposed by of Mungle et al. (2013) and Balouka and Cohen (2021), possesses the following key features:

**Variables**

| | |
|---|---|
| $x_{im} = \begin{cases} 1 \\ 0 \end{cases}$ | If Activity $i$ is implemented in Mode $m$, $x_{im}$ equals 1; otherwise, 0. |
| $y_{it} = \begin{cases} 1 \\ 0 \end{cases}$ | If Activity $i$ is completed at Time $t$, $y_{it}$ equals 1; otherwise, 0. |
| $d'_{im}$ | The real-time of implementing Activity $i$ in Mode m |
| $p_j$ | The $j^{th}$ payment rate |
| $p_J$ | The last payment rate |
| $E_i$ | Start time of Activity $i$ |
| $T_{i_j}$ | Finish time of activity $i_j$ |
| $Q_{min}$ | Minimum quality among the qualities of the selected modes for project activities |
| $Q_{avg}$ | The average quality of the selected modes for implementing project activities |

**Indices and sets**

| | |
|---|---|
| $i$ | Index of activities ($i = 1,2, \ldots, n$) |
| $j$ | Index of payment activities ($j = 1,2, \ldots, J$) |
| $m$ | Index of activity mode ($m = 1,2, \ldots, M$) |
| $t$ | Index of time ($t = 1,2, \ldots, T$) |
| $I$ | Set of the project's activities |
| $R$ | Set of the required renewable resources |
| $M$ | Set of possible implementation modes of activities |
| $succ_i$ | Set of succedent activities of Activity $i$ |
| $S_{pay}$ | Set of activities where payments are made |

**Parameters**

| | |
|---|---|
| $i_j$ | Payment j at the end of Activity $i$ |
| $c_{im}$ | The direct cost of implementing Activity $i$ under Mode $m$ |
| $R_{im}$ | Cost slope of implementing Activity $i$ under mode $m$ |
| $dem_{imr}$ | The demand of Activity $i$ for Renewable Resource r under Mode $m$ |
| $k_x$ | Interest rate |
| $D_{im}$ | Typical implementation time of Activity $i$ under Mode m according to the CPM method |
| $d_{im}$ | Implementation time of Activity $i$ with the maximally possible compression |
| $H$ | Overhead cost |
| $a_r$ | The accessible magnitude of Renewable Resource $r$ |
| $\gamma$ | Prepayment ratio |
| $\theta$ | Compensation ratio |
| $D$ | Project deadline |
| $U$ | Project price |
| $V_i$ | Value earned by Activity $i$ |
| ICA | Initial capital |
| $\alpha$ | The deviation between the minimum and average quality of the whole project |
| $q_{im}$ | The quality of Executive Mode $m$ for Activity $i$ |

Model 1: Time Cost tradeoff Problem

$$\min \sum_{m=1}^{M} \sum_{i=1}^{n} \left( \left[ \frac{C_{im} \cdot x_{im}}{(1+k_x)^{E_i + d'_{im}}} \right] + \left[ \frac{R_{im}(D_{im} \cdot x_{im} - d'_{im})}{(1+k_x)^{E_i + d'_{im}}} \right] \right) + \frac{H(E_n - E_1)}{(1+k_x)^{E_n}} \quad (4)$$

$$\min \sum_{t=EF_n}^{LF_n} y_{nt} \cdot t \quad (5)$$

s.t.

$$\sum_{m=1}^{M} x_{im} = 1 \qquad \forall i = 1,2,\dots,n \qquad (6)$$

$$\sum_{t=0}^{T} y_{it}.t = 1 \qquad \forall i = 1,2,\dots,n \qquad (7)$$

$$E_i = \sum_{t}^{T} y_{it}.t - \sum_{m=1}^{M} d'_{im} \qquad \forall i = 1,2,\dots,n \qquad (8)$$

$$P_j = (\theta - \gamma)\left[\sum_{i=1}^{n}\left(V_i \sum_{t=0}^{T_{ij}} y_{it}\right) - \sum_{i=1}^{n}\left(V_i \sum_{t=0}^{T_{ij-1}} y_{it}\right)\right] \qquad \forall j = 1,2,\dots,J-1 \qquad (9)$$

$$P_J = U - \left(\gamma U + \sum_{j=1}^{J-1} P_j\right) \qquad (10)$$

$$\sum_{i=1}^{n}\sum_{m=1}^{M} dem_{imr}.x_{im} \le a_r \qquad \forall r = 1,2,\dots R \qquad (11)$$

$$d_{im}.x_{im} \le d'_{im} \le D_{im}.x_{im} \qquad \forall i, m = 1,2,\dots,n \qquad (12)$$

$$\sum_{t=EF_n}^{LF_n} y_{nt}.t \le D \qquad (13)$$

$$\sum_{t=EF_i}^{LF_i} y_{it}.t \le \sum_{t=EF_h}^{LF_h} y_{ht}.t - \sum_{m=1}^{M} d_{hm}' \qquad \begin{matrix} h \in succ_i \\ \forall i = 1,2,\dots,n \end{matrix} \qquad (14)$$

$$\sum_{m=1}^{M}\sum_{i=1}^{n}\left(\left[\frac{C_{im}.x_{im}}{(1+k_x)^{E_i+d'_{im}}}\right] + \left[\frac{R_{im}(D_{im}.x_{im} - d'_{im})}{(1+k_x)^{E_i+d'_{im}}}\right]\right) + \frac{H(E_n - E_1)}{(1+k_x)^{E_n}}$$
$$\le ICA + \left[\gamma U + \sum_{ij \in S_{pay}} \frac{P_j \sum_{t=0}^{T} y_{ijt}}{(1+k_x)^{E_{ij}}}\right] \qquad (15)$$

$$y_{it}, x_{im} \in \{0,1\} \qquad (16)$$

In the model above, the first objective function (4) indicates the minimization of the net present value of the total costs of the project. The second objective function (5) depicts the minimization of the project implementation period. Constraint (6) implicates that every activity should be performed only in a single mode. Constraint (7) determines that every activity is initiated continuously and terminates at a specific time. Constraint (8) specifies the onset of every activity. Constraint (9) determines the $j^{th}$ payment rate. The number of the $j^{th}$ activity is calculated by the $i_j = \{i: y_{iT} = 1, \ T = min\left[T:T \geq j\left(\frac{D}{J}\right)\right]\}$ formula. With this formula, we find out in which activities we have payment action to calculate $T_{i_j} = E_{i_j} + \sum_{m=1}^{M} d'_{i_j m}$. Constraint (10) explains the rate of the last payment, which should be made at the end of the project. Constraint (11) controls resource limitations. Constraint (12) indicates that the practical duration of every activity implementation falls between the normal and compressed times. Constraint (13) denotes that the project implementation should not exceed a certain duration. Constraint (14) explains the interrelationships of the activities. Constraint (15) claims that the project implementation should be cost-effective, and constraint (16) reflects the systemic limitation of the model.

## 3.2. Productivity Enhancement

A prominent subject in project management, particularly within the context of construction projects, revolves around defining the minimum quality standards essential for a project and ensuring that these standards are met by all stakeholders and project managers. The importance of project output quality is underscored by the fact that any shortcomings in this regard can raise questions about the competence of the entire project execution team and lead to customer dissatisfaction. Consequently, one viable strategy for enhancing productivity is to elevate the overall quality of the project. Achieving this involves measuring the quality of project activities by establishing specific qualitative metrics for each task. Subsequently, following the fundamental definition of productivity (output/input ratio), a suitable productivity index can be derived by dividing the quality score obtained from project implementation by the total project costs, including both direct and indirect expenses. The ultimate objective is to improve and optimize this productivity index.

Quality, though a subjective concept, can be made measurable by breaking down the overarching quality objectives into primary properties and criteria that relate to the project's activities. Identifying these quality indices for each project activity is crucial in assessing the quality of construction work. These criteria are then assigned weights using methodologies outlined in existing literature, and the productivity index is computed. This research utilizes the DANP (Decision-Making Trial and Evaluation Laboratory) method to assign weights to qualitative measures.

## 3.3. The DANP (DEMATEL-based ANP) Method

DANP, or the DEMATEL-based ANP method, is a decision-making technique that operates with multiple criteria and sub-criteria. It constructs an Analytic Network Process (ANP) supermatrix and derives the weights for these criteria and sub-criteria through the use of the DEMATEL (Decision-Making Trial and Evaluation Laboratory) communications matrix. Essentially, DANP combines aspects of both the ANP and DEMATEL methods to form a comprehensive approach. Figure 3 illustrates the key phases of the DANP method. Within the DANP methodology, the network's structure and the weighting of dimensions are established with the assistance of the DEMATEL method. The total impact matrix generated by DEMATEL is employed to create a supermatrix without weights for ANP analysis. For more detailed information about the DANP method, you can refer to the work by Chen et al. (2014).

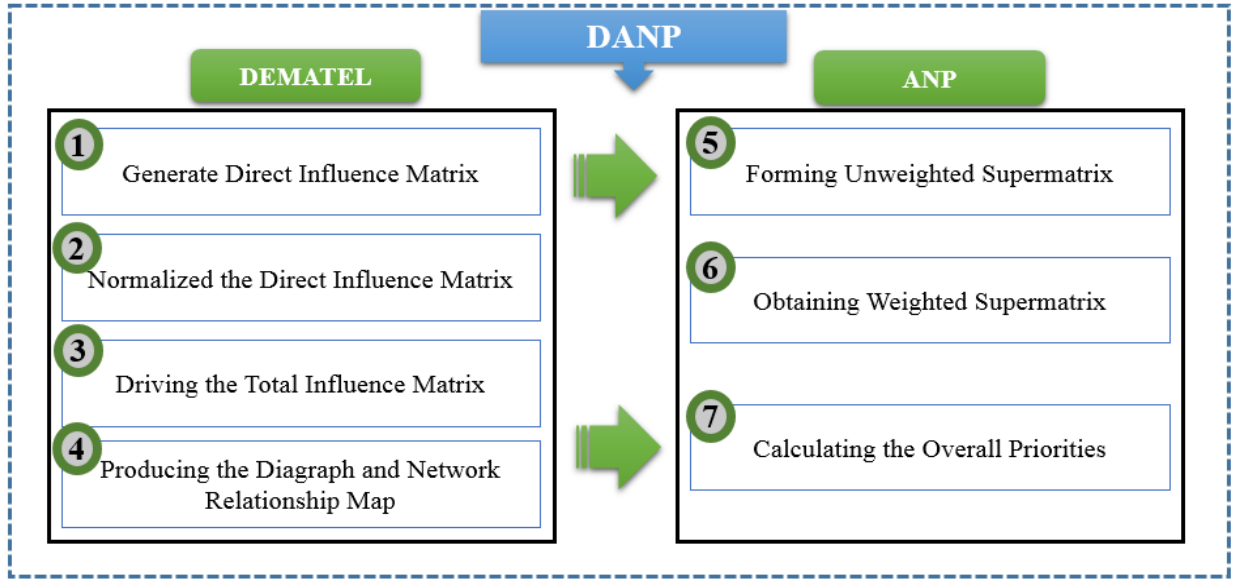

**Figure 3. DANP Method**

Now, after the calculation of the weight of the qualitative criteria and addition of the productivity improvement objective to the research problem, the three-objective non-linear mathematical model of the research takes the form below:

Model 2: Final model

$$\min \sum_{m=1}^{M} \sum_{i=1}^{n} \left( \left[ \frac{C_{im}.x_{im}}{(1+k_x)^{E_i+d'_{im}}} \right] + \left[ \frac{R_{im}(D_{im}.x_{im} - d'_{im})}{(1+k_x)^{E_i+d'_{im}}} \right] \right) + \frac{H(E_n - E_1)}{(1+k_x)^{E_n}} \quad (17)$$

$$\max \frac{\alpha * Q_{\min} + (1-\alpha)Q_{\text{avg}}}{\sum_{m=1}^{M} \sum_{i=1}^{n} \left( \left[ \frac{C_{im}.x_{im}}{(1+k_x)^{E_i+d'_{im}}} \right] + \left[ \frac{R_{im}(D_{im}.x_{im} - d'_{im})}{(1+k_x)^{E_i+d'_{im}}} \right] \right) + \frac{H(E_n - E_1)}{(1+k_x)^{E_n}}} \quad (18)$$

$$\min \sum_{t=EF_n}^{LF_n} y_{nt}.t \quad (19)$$

$s.t.$

$$Q_{\min} = \min\{q_{im}: x_{im} = 1\} \quad \begin{array}{l} 1 \leq i \leq n \\ 1 \leq m \leq M \end{array} \quad (20)$$

$$Q_{min} \leq q_{im} * x_{im} + M(1 - x_{im}) \qquad \begin{aligned} 1 &\leq i \leq n \\ 1 &\leq m \leq M \end{aligned} \qquad (21)$$

$$Q_{avg} = \frac{\sum_{m=1}^{M} \sum_{i=1}^{n}(q_{im} * x_{im})}{n} \qquad \begin{aligned} 1 &\leq i \leq n \\ 1 &\leq m \leq M \end{aligned} \qquad (22)$$

Constraints $(6) - (16)$

In the model above, the first objective function (17) indicates the minimization of the net present value of the total costs of the project. The second objective function (18) shows the productivity improvement of the project. The third objective function (19) depicts the minimization of the project implementation period. The minimum value is selected from the modes chosen for activities (20)-(21). Then, the average quality of the executive modes of activities is calculated (22).

## 4. Proposed Multi-Objective Genetic Algorithm

Considering the inherently complex nature of proposed problem, a multi-objective genetic algorithm is developed specially tailored to tackle large-scale problems. This algorithm offers several advantages for decision-makers, primarily due to its capacity for delivering effective and rapid performance, along with the ability to produce multiple solutions in each iteration, thanks to its population-oriented approach, as detailed by Ghannadpour et al. (2014).

The MOGA (Multi-Objective Genetic Algorithm) algorithm proves particularly well-suited for addressing problems characterized by high computational intricacies, such as the Multimode Resource-Constrained Project Scheduling Problem (MRCPSP). It demonstrates noteworthy efficiency, especially in the following scenarios:

**1. Large Solution Spaces:** MOGA excels when dealing with expansive solution spaces where traditional linear and non-linear programming approaches struggle to find viable solutions.

**2. Multiple Constraints:** In situations involving numerous constraints, MOGA remains robust and reliable.

**3. Time or Resource Limitations:** It thrives even when confronted with strict time limitations or other resource constraints.

**4. Approximate Solutions:** MOGA is a valuable choice when the objective is to identify approximate solutions that meet the problem's requirements.

The subsequent step involves comparing the outcomes produced by the MOGA algorithm with those generated by NSGA-II (Non-Dominated Sorting Genetic Algorithm II) algorithm developed by Deb et al. (2002). It should be noted that the NSGA-II algorithm implemented in this paper is adapted to the characteristics of the proposed problem. For a visual representation of solution algorithm, please refer to Figure 4.

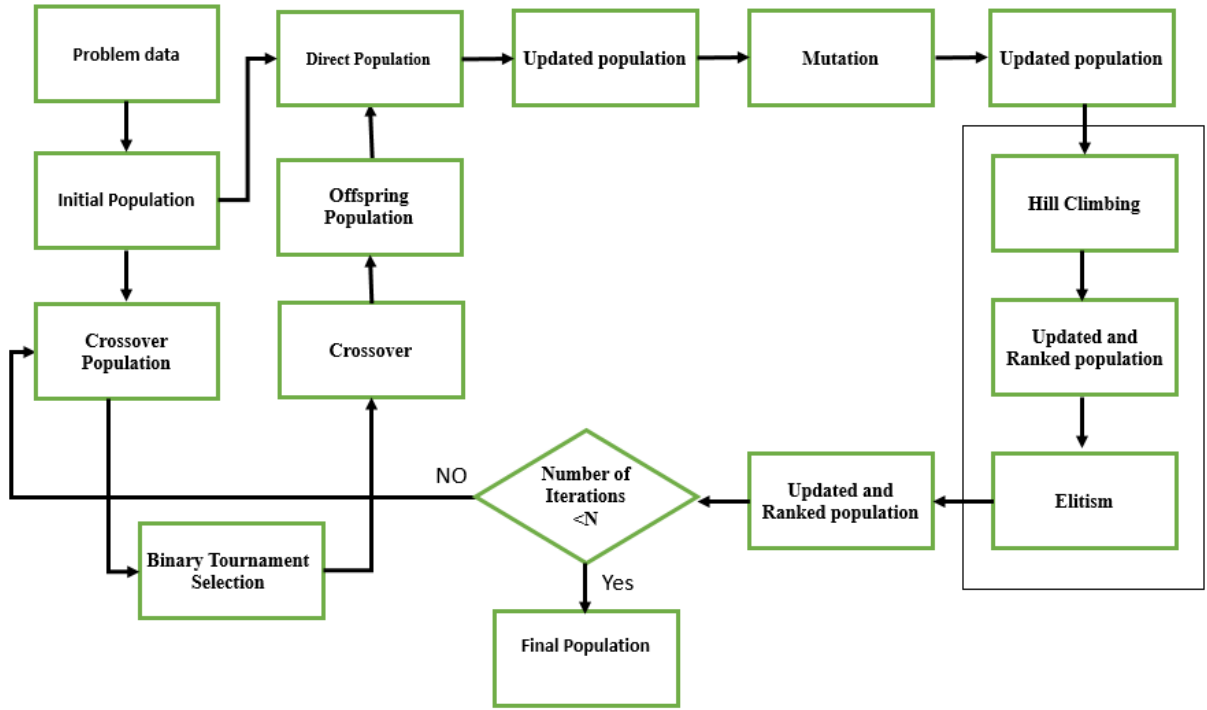

**Figure 4. Flowchart of MOGA**

## 4.1. Solution representation

In the MOGA (Multi-Objective Genetic Algorithm) algorithm employed in this research, the representation of solutions relies on a chromosome structure. Specifically, each solution is encoded using three series of integers, which serve as the genes within the chromosome. These three series convey essential information about the activities in the project.

1. Activity Sequence: The first series of integers within the chromosome represents the sequence of activities. It effectively illustrates the prerequisites and the successor relationships among activities, providing a clear depiction of their order.
2. Executive Mode: The second series of integers defines the executive mode chosen for each activity. This information indicates the specific approach or method employed in executing a particular task.
3. Activity Duration: The third series of integers indicates the time required for each activity to be completed.

For a visual depiction of this solution representation and the formation of chromosomes for the proposed problem, please refer to Figure 5. This structured representation ensures that the genetic algorithm can effectively process and optimize the solutions in line with the research objectives.

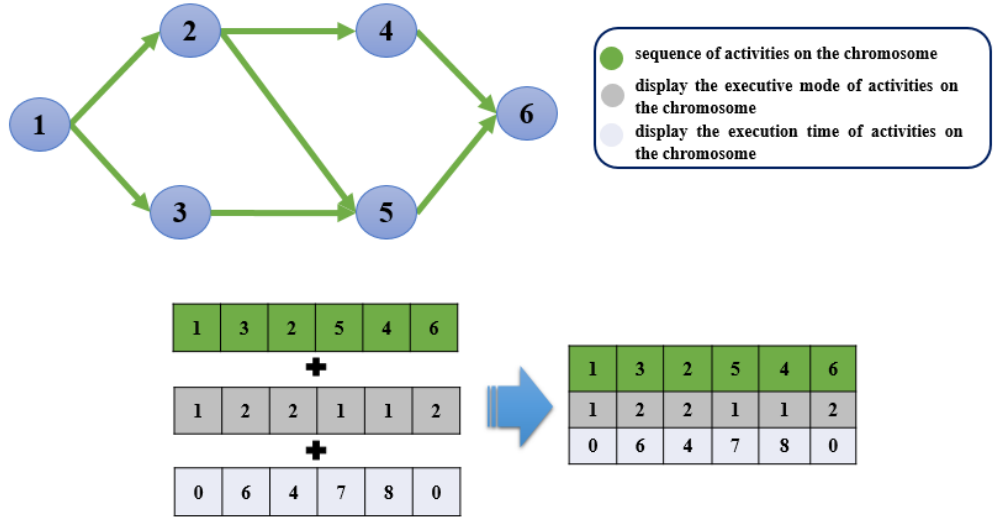

**Figure 5. Solution Representation**

## 4.2. Initial population

In addressing the proposed problem, the initial population is assembled through a random selection process. This procedure involves the generation of chromosomes, each of which is subjected to validation. If a generated chromosome adheres to the required criteria and is considered legitimate, it is incorporated into the initial population. This process iterates until the initial population is fully constituted.

For a comprehensive overview of the steps involved in the formation of the initial population, you can refer to Algorithm 1, which presents the pseudocode outlining this crucial aspect of this research. This approach ensures that the initial population is diverse and encompasses a range of potential solutions for subsequent analysis and optimization.

| Algorithm 1: Generate Random Population, Size N (pop size) for MOGA |
|---|
| 1  **for m=1:pop size** |
| 2      Valid-number = 0; |
| 3      **while valid-number** ~= 3 (Resource constraint, Time scheduling ,and  Payment and cost Constraint = 3) |
| 4          individual = Generate another random solution ; |
| 5          [valid-number, individual] = validation-function (individual,model); |
| 6      **end** |
| 7      pop{$m$, **1**} = individual.position; |
| 8      pop{$m$, **2**}  = individual.obj; |
| 9  **end** |

## 4.3. Selection method

During this phase, the selection of parent chromosomes for breeding is carried out from the population established in the previous step. The process of parent selection is governed by the tournament selection method. Here's how it works: Random Pairing: Initially, two chromosomes are randomly chosen from the population that is prepared for combining. Rank-Based Selection: Out of these two randomly selected chromosomes, the one with a lower rank within the Pareto solution set is designated as the chosen parent. If, by chance, these two chromosomes happen to be equally ranked, one is selected randomly among them. This selection process is performed twice to obtain two parents for the subsequent phase. These chosen parent chromosomes will be utilized in the breeding phase, contributing to the creation of the next generation.

### 4.4. Crossover

In the genetic algorithm, the crossover operator plays a pivotal role in producing offspring from the selected parent chromosomes. Following the establishment of the parent population, these parents are paired off, and from each pair, two offspring are generated. In this algorithm, the first string of every parent, representing the sequence of activities, is directly and immutably transferred to the respective child. However, the other two strings, one related to the representation of the executive mode of activities and the other associated with the presentation of activity implementation times within the chromosome, are shuffled between the two parents. This approach, akin to other crossover methods, yields two offspring from each pairing and introduces them into the evolving population. Figure 6 provides a visual depiction of the schematic representation of the crossover operator within the research algorithm. This operator is instrumental in diversifying the genetic makeup of the population and exploring potential solutions.

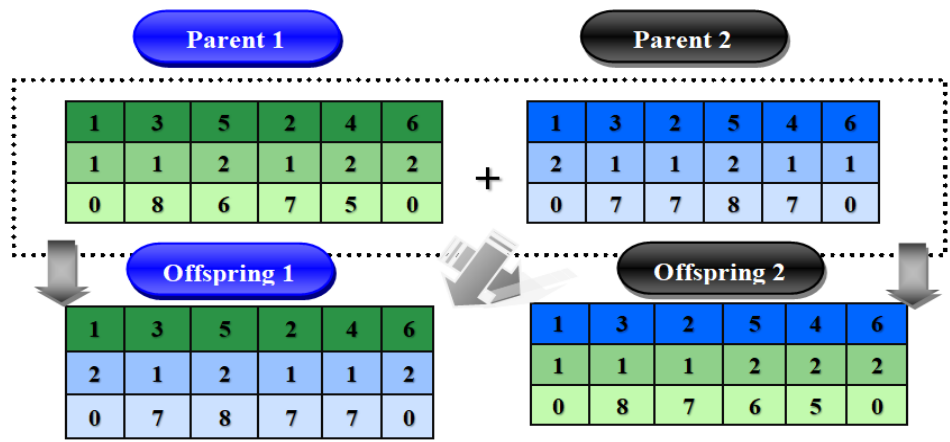

**Figure 6. Crossover operator of the MOGA algorithm**

### 4.5. Mutation

In each iteration, a portion of the new population undergoes mutation to introduce diversity among potential solutions. The mutation operator injects randomness into the genetic algorithm's search process, which is a crucial mechanism for preventing the algorithm from getting stuck in local optima. The mutation operator utilized in this research is characterized by the following features:

- **Modification of Activity Sequence:** To initiate the mutation process, two random points are selected within the chromosome representing the activity sequence (found in the first string).

- **Replacement of Executive Modes:** Next, the executive modes associated with these two selected activities, recorded in the second string of the chromosome, are replaced with new values.
- **Time Interval Adjustment:** Finally, the executive modes for these activities, located in the third string of the chromosome, are randomly reassigned within a permissible time interval.

Figure 7 provides a visual representation of how this mutation algorithm operates on a chromosome, demonstrating its role in introducing controlled variability and potentially yielding novel solutions that enhance the population's genetic diversity.

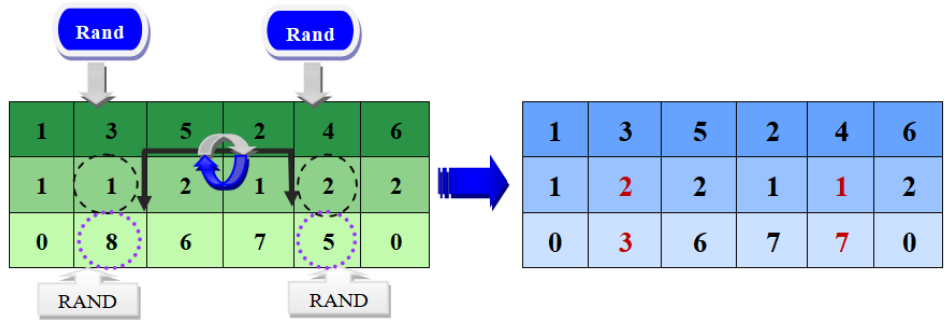

Figure 7. Mutation operator in the MOGA algorithm

## 4.6. Improvement

Following the creation of the new population, this phase focuses on enhancing the generated chromosomes. It incorporates two crucial concepts: Hill-Climbing and Elitism, both contributing to population improvement.

### 4.6.1. Hill-Climbing

Hill-Climbing is a mechanism within the improvement phase where a selected chromosome undergoes iterative refinement and is subsequently substituted for the corresponding chromosome in the population. The Hill-Climbing algorithm follows a specific set of operations to enhance the quality of the chosen chromosome.

1. **Initiation:** The process begins with the first activity in the first string of the algorithm, where the executive mode of this activity is modified.
2. **Chromosome Generation:** New chromosomes are generated, commencing from the first executable time of the activity according to its newly adjusted mode. If the generated chromosome complies with the required validity criteria, it is added to the new population. The primary chromosome remains as a member of this new population.
3. **Ranking and Comparison:** Following the ranking of the new population, if the generated chromosome secures a superior rank compared to the primary chromosome, the algorithm concludes. The improved chromosome takes the place of the initially selected one in the genetic algorithm population.

4. **Iterative Process:** If the generated chromosome does not achieve an improved rank (after assessing all activities), the same primary chromosome is reinstated in the population. The Hill-Climbing algorithm then proceeds to the next activity and repeats the entire process.

Algorithm 2 outlines the step-by-step procedure of the Hill-Climbing method, while Figure 8 offers a visual representation of how this section operates within the algorithm. Hill-Climbing aims to iteratively refine selected chromosomes, enhancing their quality and potentially leading to better solutions within the population.

| | **Algorithm 2: Hill-Climbing for MOGA** |
|---|---|
| 1 | **Choose a solution ($P$) for Hill Climbing** |
| 2 | **for $m = 2 : pop\ size - 1$**      (first and last activity have 0 duration) |
| 3 | Change the mode of activity m; $counter = 0$ |
| 4 | **for t = Crash time of activity "$m$" at new mode : Normal time of activity "$m$" at new mode** |
| 5 | Time of activity "$m$" at solution $P = t$;   (Change the time of activity m and generate "$new\ P$") |
| 6 | [$validnumber$, $new\ P$] = validation-function ($new\ P$, Model); |
| 7 | **if $validnumber == 3$** |
| 8 | pop{$cunter, 1$} = p1{$1, 1$}; |
| 9 | pop{$cunter, 2$} = p1{$cunter, 2$}; |
| 10 | $counter = counter + 1$ ; |
| 11 | **end** |
| 12 | **end** |
| 13 | Generate new pop for all new valid solution based on $P$; |
| 14 | Ranked this new population based on Pareto solutions ; |
| 15 | **If any improved is occurred** |
| 16 | $P$ substitute by improved solution ; |
| 17 | **else** |
| 18 | Continue; |
| 19 | **end** |
| 20 | **end** |
| 21 | **If nothing is changed** |
| 22 | $P$ is returned to the population ; |
| 23 | **end** |

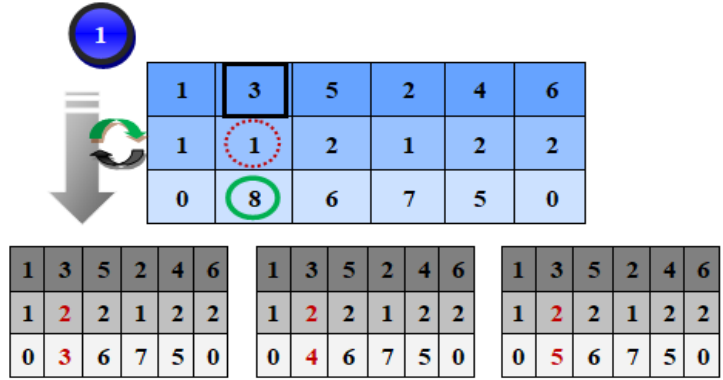

**Figure 8. Schematics of Hill Climbing operator in MOGA algorithm**

### 4.6.2. Elitism

The Elitism operator serves as a mechanism to safeguard the best solutions during each recursion of the genetic algorithm. Its operation involves retaining a portion of the top-performing solutions at the conclusion of every iteration. In practical terms, this means that a certain percentage of the best solutions from the previous population is retained, while an equivalent number of the least favorable solutions from the new population are removed. This process is supported by the work of Ghannadpour and Zandiyeh (2020) and Moosavi Heris et al. (2022).

### 4.7. Control operators

To prevent repetitious and unlegislated solutions, this research considers a number of control operators before and after the crossover, mutation, and improvement operators and disregards and eliminates the recurrent member if there exists any.

## 5. Numerical results

### 5.1. Parameter tuning

In this section, the Taguchi method is employed to fine-tune the parameters of the MOGA algorithm, capitalizing on its advantages including the reduction in required tests, cost and time savings, the capacity to explore interactive effects and conduct parallel tests, and the ability to predict optimal solutions. Quality assessment is based on deviations from desired values, and the Analysis of Mean (ANOM) approach optimizes results obtained from a single round of testing. Mean values for each factor at different levels are estimated using Formula (23), leading to the determination of optimal levels and combinations through solution tables and diagrams. This systematic parameter tuning process enhances the MOGA algorithm's effectiveness in solving complex optimization problems.

$$(M)_{Factor=I}^{Level=i} = \frac{1}{n}\sum [(f)_{Factor=I}^{Level=i}] \qquad (23)$$

Parameter tuning in the MOGA algorithm encompasses six parameters evaluated at five different levels, as outlined in Table 2. To facilitate this optimization, the standard Taguchi table, specifically the L25 test with 25 experiments, is employed. Once these Taguchi-recommended tests are conducted, the next crucial step involves identifying and analyzing the parameters that significantly impact the algorithm's performance quality. This analysis aims to pinpoint the optimal parameter values and extract the best solutions derived from the Taguchi tests, ultimately contributing to the fine-tuning and enhancement of the algorithm's effectiveness.

**Table 2. Parameters and different levels of Taguchi test of MOGA algorithm**

| Levels | Elitism | Hill Climbing | Mutation | Crossover | Iterations | Number of population |
|--------|---------|---------------|----------|-----------|------------|----------------------|
| 1 | 0.05 | 0.2 | 0.2 | 0.2 | 500 | 20 |
| 2 | 0.1 | 0.4 | 0.4 | 0.4 | 700 | 40 |
| 3 | 0.15 | 0.5 | 0.5 | 0.5 | 1100 | 60 |
| 4 | 0.2 | 0.6 | 0.6 | 0.6 | 1500 | 80 |
| 5 | 0.25 | 0.8 | 0.8 | 0.8 | 2000 | 100 |

After designing and conducting the 25 Taguchi tests and calculating the relevant criteria, the test results are obtained. Table 3 and Figure 9 present the outcomes of the parameter tuning process applied to the MOGA algorithm. These results provide insights into the optimized parameter values and their impact on the algorithm's performance, further refining its effectiveness.

**Table 3. The optimum level of parameters for the algorithm**

| | Elitism | Hill Climbing | Mutation | Crossover | Iterations | Number of population |
|---|---------|---------------|----------|-----------|------------|----------------------|
| Optimum Level | 0.05 | 0.8 | 0.6 | 0.8 | 2000 | 100 |

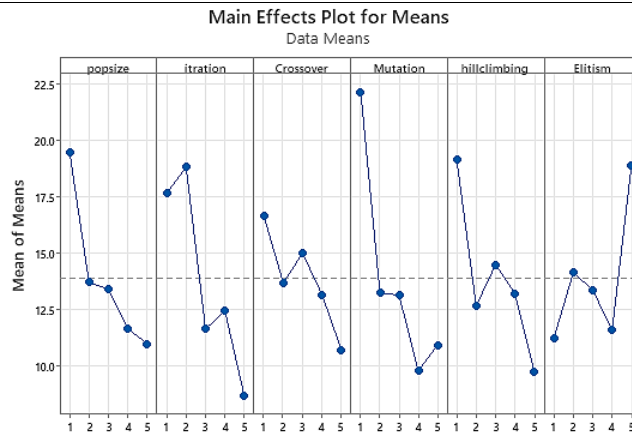

**Figure 9. Mean graph for MOGA algorithm at different factor levels**

The algorithm validation process consists of two phases. Firstly, small-size problems are solved using a precise method with the BARONS solver in GAMS software, and their results are compared with those obtained from the MOGA algorithm. Secondly, for larger-scale problems, the MOGA and NSGA-II algorithms are employed, and their outcomes are contrasted. This validation ensures the algorithm's reliability and effectiveness across varying problem sizes and complexities.

## 5.2. Small-Scale Problem Validation

To validate the MOGA algorithm, small-scale problems are considered in this section. These problems fall under the category of MMRCPSP, which is known to be unsolvable by accurate methods, especially for larger sizes and within limited time frames. To validate the MOGA algorithm's effectiveness, these small-scale problems are defined and solved using the precise ε-constraint method via the BARON solver within the GAMS software. Subsequently, the results obtained from this accurate method are compared with the solutions generated by the proposed MOGA algorithm. In total, five different problems are examined, and their respective results are presented in Table 4. The comparison reveals that both the accurate method and the MOGA algorithm produce nearly identical values for the objective functions. The key distinction lies in the fact that the genetic algorithm's computational time is considerably shorter than that of the accurate method in GAMS software. It's important to note that the accurate method reaches its limitations for problems with nine dimensions, taking over two hours to solve, which led to the software stopping.

**Table 4. GAMS and MOGA results in solving small-Scale problems**

| No. | Num. of Activities | MOGA | | | | | GAMS | | | | |
|---|---|---|---|---|---|---|---|---|---|---|---|
| | | Num. of Pareto | NPV of Project Cost | Time | Productivity Index | CPU time (S) | Num. of Pareto | NPV of Project Cost | Time | Productivity Index | CPU time (S) |
| 1 | 6 | 18 | 475.12 | 10 | 0.22 | 312 | 8 | 475.12 | 10 | 0.22 | 2235 |
| 2 | 6 | 19 | 366.7 | 9 | 0.28 | 325 | 10 | 366.7 | 9 | 0.28 | 2385 |
| 3 | 6 | 19 | 388.36 | 10 | 0.25 | 342 | 9 | 388.36 | 10 | 0.25 | 2555 |
| 4 | 7 | 25 | 392.26 | 14 | 0.25 | 389 | 12 | 392.26 | 14 | 0.25 | 4555 |
| 5 | 9 | 17 | 646.2 | 9 | 0.14 | 455 | - | - | - | - | - |

Problem 4 is subjected to a more detailed analysis, as illustrated in Figure 10 and Figure 11. In the optimal solution related to the project completion time, this project concludes in a duration of 14 time units. However, the net present value of the project costs exhibits an opposite trend, increasing. Additionally, the overall productivity rate decreases. Conversely, in the optimal solution associated with the net present value of the project costs, the objective function attains a value of 392.9. However, the project completion time departs from its optimal state, extending to 25 time units. Remarkably, the net present value of the project costs aligns with productivity improvement, also witnessing an increase. The overall productivity enhancement of the project remains consistent across the third objective function. This analysis provides insights into the trade-offs and interactions between these critical project parameters, shedding light on the complex dynamics of problem 4.

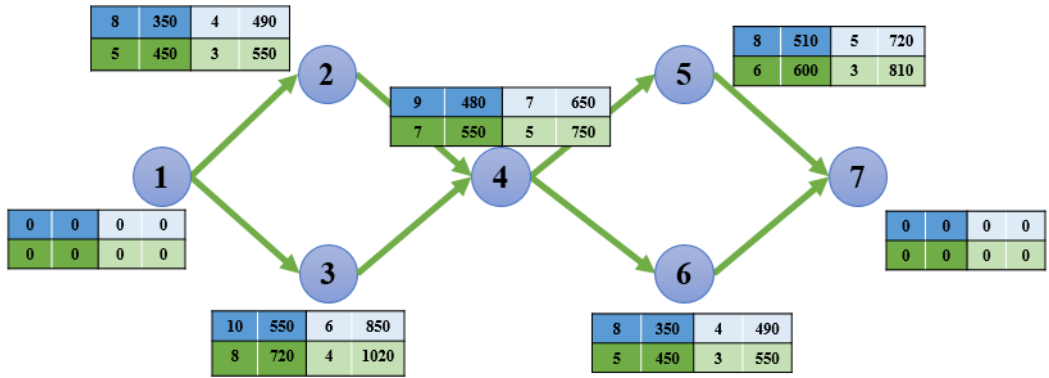

**Figure 10. Analysis of problem-solving considering 7 dimensions**

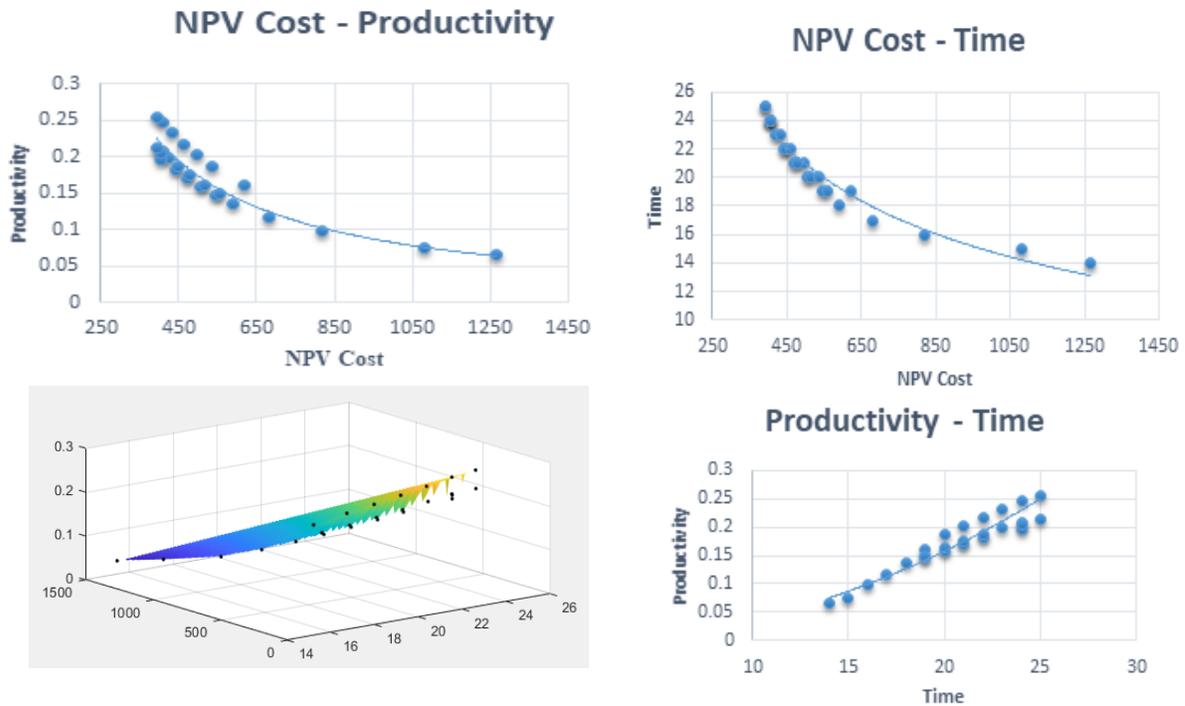

**Figure 11. Pareto frontier developed by MOGA algorithm**

## 5.3. Large-Scale problems validation

This research conducts an in-depth analysis of problems with large dimensions to evaluate the performance of the MOGA algorithm. Given the inherent complexity of the research problem, a set of 23 problems is meticulously designed, varying in size, with relevant data integrated into the project Graph following the mathematical model used in this study. To facilitate a robust performance comparison, these problems are also solved using the NSGA-II algorithm, and the results are systematically compared with those generated by MOGA. NSGA-II, functioning on the principles of a genetic algorithm, serves as a standard benchmark for assessing the performance of other multi-objective meta-heuristic algorithms. To ensure a fair comparison, parameter tuning is applied to NSGA-II, ensuring that both algorithms operate under optimal conditions. The comparison between the two algorithms is based on the following criteria:

### 5.3.1. Best Solution (BS)

This criterion encompasses three aspects, including the solution with the best net present value of costs, the solution with the best project implementation time, and the solution with the best productivity. For the first two objective functions, which involve time and cost, smaller values indicate superior algorithm performance as these objectives are of minimization type. Conversely, the third objective function, productivity, is of maximization type, and larger values signify better algorithm performance in this context. This multifaceted evaluation provides a holistic view of the algorithm's effectiveness

### 5.3.2. The Mean Ideal Distance (MID) criterion

The Mean Ideal Distance (MID) criterion assesses the proximity of Pareto solutions to the optimal solution. Smaller values of this index indicate that the solutions are closer to the optimum solution, signifying the algorithm's effectiveness. The MID is calculated using Formula (24), where *n* represents the number of Pareto solutions. For further details on the MID criterion, readers can refer to Fathollahi-Fard et al. (2020). This criterion provides valuable insights into the algorithm's ability to approximate optimal solutions.

$$\text{MID} = \frac{\sum_{i=1}^{n}\sqrt{\left(\frac{f_{1i}-f_{1\,best}}{f_{1\,max}-f_{1\,min}}\right)^2 + \left(\frac{f_{2i}-f_{2\,best}}{f_{2\,max}-f_{2\,min}}\right)^2 + \left(\frac{f_{3i}-f_{3\,best}}{f_{3\,max}-f_{3\,min}}\right)^2}}{n} \tag{24}$$

### 5.3.3. Diversity Measure (DM)

This index measures the diversity of Pareto solutions obtained by the algorithm. The larger values of this index depict the proximity of the solutions to the optimum solution and thus the further utility of the algorithm. Refer Majumdar et al. (2022) for more information about the DM

### 5.3.4. Number of Pareto Solutions (NPS) criterion

This index represents the number of Pareto solutions obtained in the Pareto frontier. The larger values of this index indicate the proximity of the solutions to the optimum solution and thus the further utility of the algorithm (Sadeghi et al., 2021).

### 5.3.5. Quality Measure (QM)

For the estimation of this criterion, the Pareto solutions of all algorithms are mixed, and a new Pareto set is developed. Then, the number of the Pareto solutions of every algorithm, belonging to the combinatorial set of the Pareto solutions, is calculated. The algorithm with a higher value of this criterion reveals better performance (Shang et al., 2021).

### 5.3.6. Generational distance (GD)

The GD indicator captures the average distance between each element of a Pareto front approximation and its closest neighbor in a discrete representation of the Pareto front. This indicator is given by Formula (25).

$$GD(Y_N : Y_P) = \frac{1}{|Y_N|}\left(\sum_{y_1 \in Y_N} min_{y_2 \in Y_P} \|y^1 - y^2\|^p\right)^{\frac{1}{p}} \tag{25}$$

where $|Y_N|$ is the number of points in a Pareto front approximation and $Y_P \subseteq y_p$ a discrete representation of the Pareto front. For all these indicators, a lower value is considered to be better (Audet et al., 2021).

### 5.3.7. Maximum Pareto front error (MPFE):

This indicator is another measure that evaluates the distance between a discrete representation of the Pareto front and the Pareto front approximation obtained by a given algorithm. It corresponds to the largest minimal distance between elements of the Pareto front approximation and their closest neighbors belonging to the Pareto front. This indicator is to be minimized. It is expressed with Formula (26) (generally, $p = 2$) (Audet et al., 2021):

$$MPFE(Y_N:Y_P) = max_{y^2 \in Y_N} \left( min_{y^2 \in Y_p} \sum_{i=1}^{m} |y_1^1 - y_1^2|^p \right)^{\frac{1}{p}} \quad (26)$$

### 5.3.8. Spacing (SP):

The SP indicator captures the variation of the distance between elements of a Pareto front approximation. A lower value is considered to be better (Audet et al., 2021). This indicator is computed with Formula (27).

$$SP(Y_N) = \sqrt{\frac{1}{|Y_N|-1} \sum_{j=1}^{|Y_N|} (\bar{d} - d^1(y^j, Y_N \setminus \{y^j\}))^2} \quad (27)$$

where $d^1(y^j, Y_N \setminus \{y^j\}) = min_{y \in Y_N \setminus \{y^j\}} \|y - y^j\|_1$ is the $l_1$ distance of $y^j \in Y_N$ to the set $Y_N \setminus \{y^j\}$ and $\bar{d}$ is the mean of all $d^1(y^j, Y_N \setminus \{y^j\})$ for $j = 1, 2, ..., |Y_N|$.

### 5.3.9. Hole relative size (HRS):

This indicator identifies the largest hole in a Pareto front approximation for a bi-objective problem. It is given by Formula (28).

$$HRS(Y_N) = (1/\bar{d}) \, max_{j=1.2...|Y_N|-1} d^j \quad (28)$$

where $Y_N$ is a Pareto front approximation whose elements are sorted in ascendant order according to the first objective, $d^j = \|y^j - y^{j+1}\|_2$ is the $l_2$ distance between the two adjacent objective vectors $y^j \in Y_N$ and $y^{j+1} \in Y_N$ and $\bar{d}$ the mean of all $d^j$ for $j = 1, 2, ..., |Y_N| - 1$. A lower indicator value is desirable. As it takes into account holes in the objective space, this indicator is more adapted to continuous Pareto front approximations (Audet et al., 2021).

The introduced criteria for comparing two algorithms are calculated for the 23 problems solved by the two algorithms, and the results are presented in Table 5 and Table 6. Likewise, the difference between the algorithms in these criteria displays in Table 7. In all introduced criteria, the genetic algorithm acquires better values than NSGA-II.

Table 5. Evaluation criteria of MOGA algorithm

| MOGA | | | | | | | | | | | | | |
|---|---|---|---|---|---|---|---|---|---|---|---|---|---|
| No. | Number of Activities | Best Productivity | NPV Cost | Time | NPS | QM | DM | MID | GD | MPFE | HRS | SP | CPU time (s) |
| 1 | 6 | 0.211 | 475.12 | 10 | 18 | 0.56 | 710.01 | 3.42 | 30/6 | 95.5 | 3.2 | 1.8 | 312 |
| 2 | 6 | 0.28 | 366.1 | 9 | 19 | 0.64 | 628.43 | 2.41 | 22.64 | 185.08 | 4.08 | 3.36 | 325 |
| 3 | 6 | 0.25 | 388.35 | 10 | 19 | 0.67 | 710.47 | 2.94 | 27.41 | 176.73 | 2.4 | 2.29 | 342 |
| 4 | 9 | 0.14 | 646.2 | 9 | 17 | 0.75 | 1336.3 | 4.45 | 68.30 | 381.02 | 4.33 | 21.35 | 455 |

| No. | Number of Activities | Best productivity | NPV Cost | Time | NPS | QM | DM | MID | GD | MPFE | HRS | SP | CPU time (s) |
|---|---|---|---|---|---|---|---|---|---|---|---|---|---|
| 5 | 9 | 0.25 | 420.1 | 17 | 31 | 0.64 | 629.53 | 2.92 | 11.43 | 161.87 | 5.64 | 10.24 | 500 |
| 6 | 9 | 0.46 | 239.3 | 17 | 34 | 0.76 | 976.68 | 1.81 | 18.96 | 212.85 | 6.03 | 25.04 | 505 |
| 7 | 12 | 0.15 | 643.9 | 19 | 13 | 0.59 | 761.45 | 5.27 | 57.27 | 279.8 | 4.03 | 5.5 | 752 |
| 8 | 12 | 0.26 | 407.4 | 23 | 33 | 0.82 | 778.01 | 2.76 | 27.56 | 86.25 | 4.62 | 5.83 | 782 |
| 9 | 15 | 0.115 | 822.7 | 17 | 22 | 0.63 | 1366.3 | 5.86 | 35.27 | 225.14 | 5.96 | 14.78 | 865 |
| 10 | 15 | 0.203 | 486.2 | 19 | 20 | 0.64 | 1371.4 | 3.51 | 33.3 | 261.4 | 3.96 | 5.3 | 900 |
| 11 | 18 | 0.16 | 645.7 | 31 | 36 | 0.66 | 895.16 | 4.25 | 27.85 | 83.12 | 4.62 | 4.87 | 1450 |
| 12 | 18 | 0.27 | 407.7 | 37 | 50 | 0.63 | 812.7 | 2.9 | 26.26 | 242.8 | 5.3 | 10.9 | 1652 |
| 13 | 18 | 0.256 | 407.7 | 37 | 50 | 0.63 | 812.7 | 2.9 | 27.52 | 88.36 | 4.56 | 7.2 | 1780 |
| 14 | 21 | 0.16 | 632.7 | 41 | 52 | 0.76 | 1029.6 | 4.09 | 15.8 | 396.35 | 10.725 | 24.5 | 2562 |
| 15 | 21 | 0.14 | 731 | 29 | 31 | 0.59 | 1021.11 | 4.77 | 48.8 | 320.2 | 4.3 | 12.5 | 2620 |
| 16 | 24 | 0.15 | 686.5 | 36 | 53 | 0.56 | 768.31 | 5.11 | 35.5 | 125.2 | 7.6 | 11.5 | 3254 |
| 17 | 24 | 0.15 | 686.5 | 35 | 51 | 0.59 | 768.31 | 5.11 | 15.25 | 87.5 | 3.25 | 5.2 | 3580 |
| 18 | 27 | 0.12 | 840.7 | 24 | 32 | 0.85 | 1533.78 | 4.66 | 20.33 | 161.1 | 7.48 | 17.24 | 3560 |
| 19 | 27 | 0.25 | 398.2 | 46 | 30 | 0.78 | 888.35 | 2.87 | 27.12 | 303.6 | 8.2 | 9.18 | 3700 |
| 20 | 30 | 0.11 | 938.5 | 27 | 24 | 0.69 | 1110.88 | 7 | 31.98 | 214.5 | 3.95 | 9.93 | 3850 |
| 21 | 30 | 0.33 | 238.5 | 47 | 33 | 0.81 | 791.84 | 2.3 | 11.19 | 87.81 | 2.95 | 4.67 | 3950 |
| 22 | 40 | 0.15 | 648.9 | 54 | 36 | 0.89 | 906.07 | 4.26 | 19.68 | 245.54 | 8.25 | 10.22 | 4205 |
| 23 | 40 | 0.26 | 394.8 | 50 | 38 | 0.80 | 813.64 | 2.87 | 22.01 | 328.46 | 13.20 | 15.2 | 4250 |

Table 6. Evaluation criteria of NSGA-II algorithm

| NSGA-II | | | | | | | | | | | | | |
|---|---|---|---|---|---|---|---|---|---|---|---|---|---|
| No. | Number of Activities | Best productivity | NPV Cost | Time | NPS | QM | DM | MID | GD | MPFE | HRS | SP | CPU time (s) |
| 1 | 6 | 0.21 | 475.12 | 10 | 18 | 0.44 | 710.01 | 3.43 | 35.12 | 95.5 | 3.42 | 2.2 | 68 |
| 2 | 6 | 0.28 | 366.1 | 9 | 18 | 0.36 | 628.43 | 2.5 | 27.56 | 185.082 | 4.08 | 4.47 | 69 |
| 3 | 6 | 0.25 | 388.35 | 10 | 12 | 0.33 | 710.47 | 2.96 | 45.34 | 176.73 | 3.04 | 2.18 | 68 |
| 4 | 9 | 0.123 | 668.05 | 9 | 12 | 0.25 | 1259 | 5.36 | 70.25 | 383.23 | 4.5 | 22.12 | 85 |
| 5 | 9 | 0.222 | 424.63 | 17 | 23 | 0.36 | 629.53 | 3.08 | 20.23 | 164.8 | 6.64 | 11.24 | 90 |
| 6 | 9 | 0.39 | 259.1 | 17 | 20 | 0.23 | 906.83 | 1.87 | 41.87 | 288.86 | 6.97 | 11.28 | 91 |
| 7 | 12 | 0.135 | 654.84 | 19 | 11 | 0.41 | 750.48 | 5.98 | 68.88 | 279.8 | 3.66 | 5.9 | 120 |
| 8 | 12 | 0.223 | 412.16 | 24 | 17 | 0.18 | 660.24 | 3.3 | 27.74 | 140.28 | 5.7 | 6.12 | 135 |
| 9 | 15 | 0.1 | 913.63 | 17 | 14 | 0.37 | 1167/56 | 7.87 | 67.95 | 565.83 | 6.9 | 44.17 | 178 |
| 10 | 15 | 0.18 | 489.66 | 19 | 14 | 0.36 | 1080.95 | 3.87 | 66.81 | 299.66 | 7.3 | 5.4 | 181 |
| 11 | 18 | 0.13 | 667.25 | 31 | 14 | 0.34 | 865.97 | 5.95 | 62.23 | 120.21 | 6.2 | 5.12 | 272 |
| 12 | 18 | 0.2 | 430.95 | 38 | 19 | 0.37 | 711.99 | 3.94 | 36.26 | 263.8 | 6.3 | 11.9 | 280 |
| 13 | 18 | 0.2 | 430.95 | 36 | 19 | 0.37 | 711.99 | 3.94 | 35.52 | 135.25 | 5.56 | 8.2 | 280 |
| 14 | 21 | 0.13 | 656.25 | 42 | 14 | 0.23 | 1008.31 | 5.98 | 56.79 | 497.8 | 16.41 | 26.9 | 900 |
| 15 | 21 | 0.111 | 795.67 | 28 | 15 | 0.41 | 1146.88 | 6.3 | 67.18 | 356.1 | 4.3 | 15.7 | 625 |
| 16 | 24 | 0.112 | 764.34 | 36 | 17 | 0.44 | 707.46 | 7.75 | 48.6 | 185.2 | 7.9 | 11.8 | 956 |
| 17 | 24 | 0.112 | 764.34 | 36 | 19 | 0.41 | 707.46 | 7.75 | 32.5 | 145.89 | 6.87 | 8.3 | 957 |
| 18 | 27 | 0.1 | 956.41 | 24 | 17 | 0.15 | 1351.03 | 7.94 | 38.47 | 187.5 | 7.8 | 17.89 | 1012 |

| No. | | | | | | | | | | | | | |
|---|---|---|---|---|---|---|---|---|---|---|---|---|---|
| 19 | 27 | 0.2 | 432 | 46 | 18 | 0.22 | 739.43 | 3.51 | 27.55 | 308.45 | 5.8 | 6.5 | 1120 |
| 20 | 30 | 0.09 | 1100 | 27 | 18 | 0.31 | 958.82 | 11.44 | 33.62 | 134.5 | 3.66 | 8.66 | 1502 |
| 21 | 30 | 0.3 | 302 | 47 | 25 | 0.19 | 670.02 | 2.6 | 14.47 | 134.4 | 6.42 | 11.67 | 1602 |
| 22 | 40 | 0.12 | 713.3 | 53 | 17 | 0.11 | 807.03 | 7.56 | 28.58 | 253.76 | 9.3 | 13.88 | 1950 |
| 23 | 40 | 0.2 | 432.9 | 50 | 13 | 0.19 | 787 | 4.03 | 55.18 | 334.03 | 14.25 | 17.2 | 2015 |

**Table 7. Comparison of the Genetic and NSGA-II algorithms**

**NSGA-II VS MOGA**

| No. | Number of Activities | Difference of NPS (%) | The Difference of the Best (%) | | | Difference of QM (%) | Difference of DM (%) | Difference of MID (%) | Difference of GD (%) | Difference of MPFE (%) | Difference of HRS (%) | Difference of SP (%) |
|---|---|---|---|---|---|---|---|---|---|---|---|---|
| | | | productivity | NPV Cost | Time | | | | | | | |
| 1 | 6 | 0.00 | 0.00 | 0.00 | 0.00 | 27.27 | 0.00 | -0.26 | -18.18 | -6.43 | 0.00 | -12.8702 |
| 2 | 6 | 5.56 | 0.00 | 0.00 | 0.00 | 77.78 | 0.00 | 0.00 | -24.83 | 0.00 | 0.00 | -17.85 |
| 3 | 6 | 58.33 | 0.86 | 0.00 | 0.00 | 100.00 | 0.00 | -0.77 | 5.05 | -21.05 | 0.00 | -39.55 |
| 4 | 9 | 41.67 | 11.56 | -3.27 | 0.00 | 200.00 | 6.14 | -17.03 | -3.48 | -3.78 | -0.58 | -2.78 |
| 5 | 9 | 34.78 | 14.43 | -1.07 | 0.00 | 78.57 | -4.93 | -4.98 | -8.90 | -15.06 | -1.78 | -43.50 |
| 6 | 9 | 70.00 | 17.13 | -7.63 | 0.00 | 225.00 | 7.70 | -3.27 | 121.99 | -13.49 | -26.31 | -54.72 |
| 7 | 12 | 18.18 | 8.86 | -1.67 | 0.00 | 42.86 | 1.46 | -11.84 | -6.78 | 10.11 | 0.00 | -16.86 |
| 8 | 12 | 94.12 | 14.80 | -1.15 | -4.17 | 350.00 | 17.84 | -16.34 | -4.74 | -18.95 | -38.52 | -0.65 |
| 9 | 15 | 57.14 | 19.68 | -9.96 | 0.00 | 71.43 | 17.00 | -25.54 | -66.54 | -13.62 | -60.21 | -48.09 |
| 10 | 15 | 42.86 | 14.33 | -2.51 | 0.00 | 75.00 | 26.87 | -9.15 | -1.85 | -45.75 | -12.77 | -50.16 |
| 11 | 18 | 157.14 | 24.16 | -3.23 | 0.00 | 91.63 | 3.36 | -28.52 | -4.88 | -25.48 | -30.85 | -55.25 |
| 12 | 18 | 163.16 | 31.25 | -5.39 | -2.63 | 71.43 | 14.14 | -26.33 | -8.40 | -15.87 | -7.96 | -27.58 |
| 13 | 18 | 163.16 | 31.25 | -5.39 | 2.78 | 71.43 | 14.14 | -26.33 | -12.20 | -17.99 | -34.67 | -22.52 |
| 14 | 21 | 271.43 | 21.84 | -3.60 | -2.38 | 225.00 | 2.11 | -31.69 | -8.92 | -34.64 | -20.38 | -72.18 |
| 15 | 21 | 106.67 | 27.30 | -8.13 | 3.57 | 42.86 | -10.97 | -24.27 | -20.38 | 0.00 | -10.08 | -27.36 |
| 16 | 24 | 211.77 | 32.74 | -10.18 | 0.00 | 27.27 | 8.60 | -34.01 | -2.54 | -3.80 | -32.40 | -26.95 |
| 17 | 24 | 168.42 | 32.74 | -10.18 | -2.78 | 42.86 | 8.60 | -34.01 | -37.35 | -52.69 | -40.02 | -53.08 |
| 18 | 27 | 88.24 | 29.09 | -12.10 | 0.00 | 460.00 | 13.53 | -41.30 | -3.63 | -4.10 | -14.08 | -47.15 |
| 19 | 27 | 66.67 | 25.75 | -7.81 | 0.00 | 257.14 | 20.14 | -18.07 | 41.23 | 41.38 | -1.57 | -1.56 |
| 20 | 30 | 33.33 | 34.76 | -14.67 | 0.00 | 122.22 | 15.85 | -38.82 | 14.67 | 7.92 | 59.48 | -4.88 |
| 21 | 30 | 32.00 | 17.92 | -6.02 | 0.00 | 314.30 | 18.20 | -11.65 | -59.98 | -54.05 | -34.67 | -22.67 |
| 22 | 40 | 111.76 | 35.88 | -9.02 | 1.89 | 700.00 | 12.24 | -43.50 | -26.37 | -11.29 | -3.24 | -31.14 |
| 23 | 40 | 192.31 | 33.62 | -8.80 | 0.00 | 312.50 | 3.38 | -28.75 | -11.63 | -7.37 | -1.67 | -60.11 |
| **Average** | | **95.17** | **20.87** | **-5.73** | **-0.16** | **173.33** | **8.50** | **-20.72** | **-6.46** | **-13.30** | **-13.58** | **-32.15** |

Table 7 provides a detailed comparison of the proposed genetic algorithm and NSGA-II across various evaluation criteria. The genetic algorithm demonstrates superior performance, with an average 20.87% improvement in the total productivity objective function compared to NSGA-II. Additionally, the net present value of project costs in the genetic algorithm is 5.73% higher than NSGA-II, while the project completion time objective function in the genetic algorithm is only 0.16% different from NSGA-II on average. In terms of diversity, the genetic algorithm outperforms NSGA-II, exhibiting an 8.5% higher diversity across the 23 problems. The MID values further emphasize the genetic algorithm's utility, with a 20.72% reduction compared to NSGA-II. Moreover, the genetic algorithm excels with a 173.33% difference in the Quality Measure and a 95.17% advantage in the NPS index over NSGA-II. Figure 12 presents a radar graph illustrating the performance differences between the two algorithms. This graph is generated by computing the average indices for each problem and normalizing their values, with proximity to one indicating algorithm outperformance. Figure 13 to 15 provide a visual comparison of the objective function values obtained by the two algorithms.

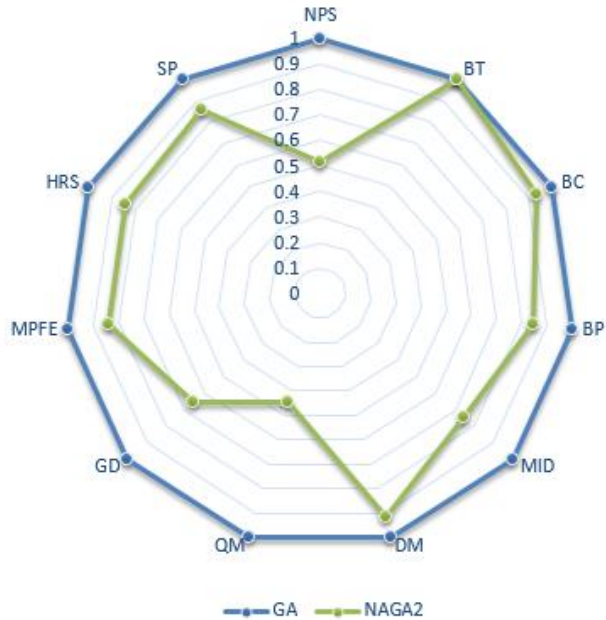

**Figure 12. Radar graph comparing two algorithms according to the criteria**

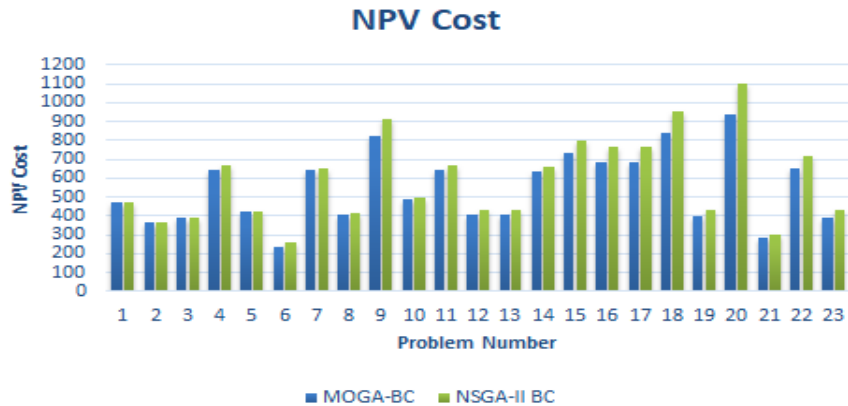

**Figure 13. Comparing the values of the net present value of costs objective function**

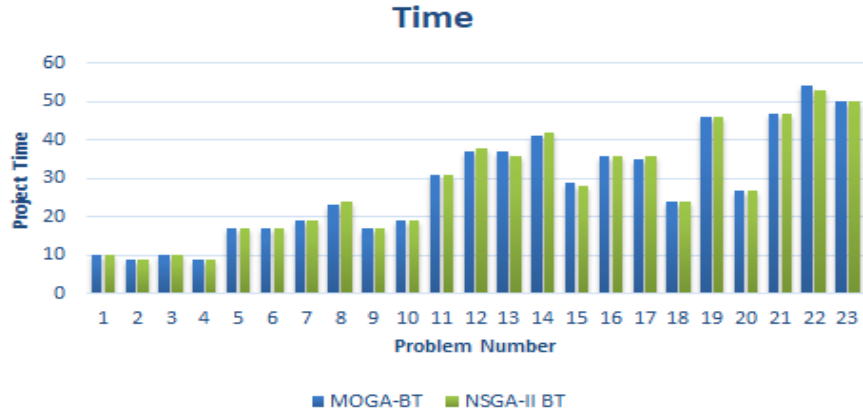

**Figure 14. Comparing the values of the project completion time objective function**

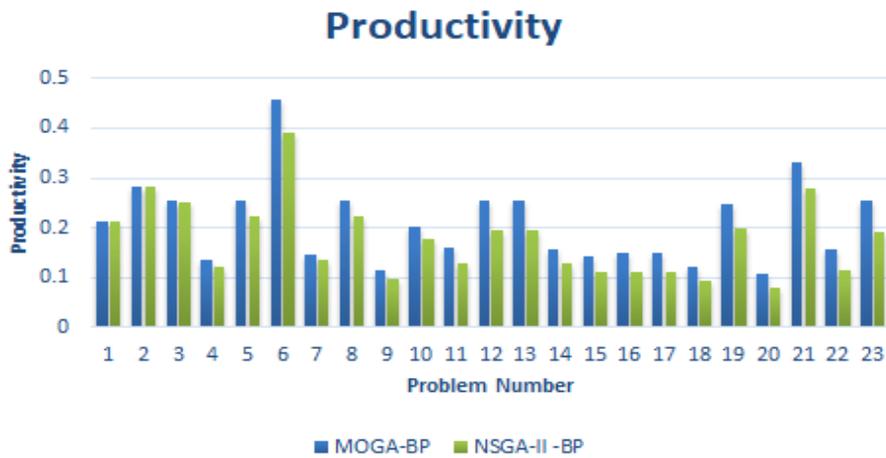

**Figure 15. . Comparing the values of the project productivity objective function**

In general, both algorithms exhibit relatively similar performance in terms of the time objective function and diversity measure. However, the genetic algorithm surpasses NSGA-II in several other criteria. Figure 16 visually depicts the Pareto frontier developed for problems featuring 27, 30, and 40 dimensions. The input data for these problems highlight that as project costs are compressed with increased expenditure, the quality of activities diminishes. This is reflected in the clear trade-offs between the cost and time functions, as well as between time and productivity functions, evident in the constructed Pareto frontiers. Specifically, when seeking to minimize the net present value of costs, the time objective function tends to increase, and vice versa. Conversely, there is no inherent trade-off between the cost and productivity functions; striving to reduce the net present value of costs tends to result in an increase in the productivity objective function and vice versa. These insights provide valuable context for understanding the relationship between project variables and optimization objectives.

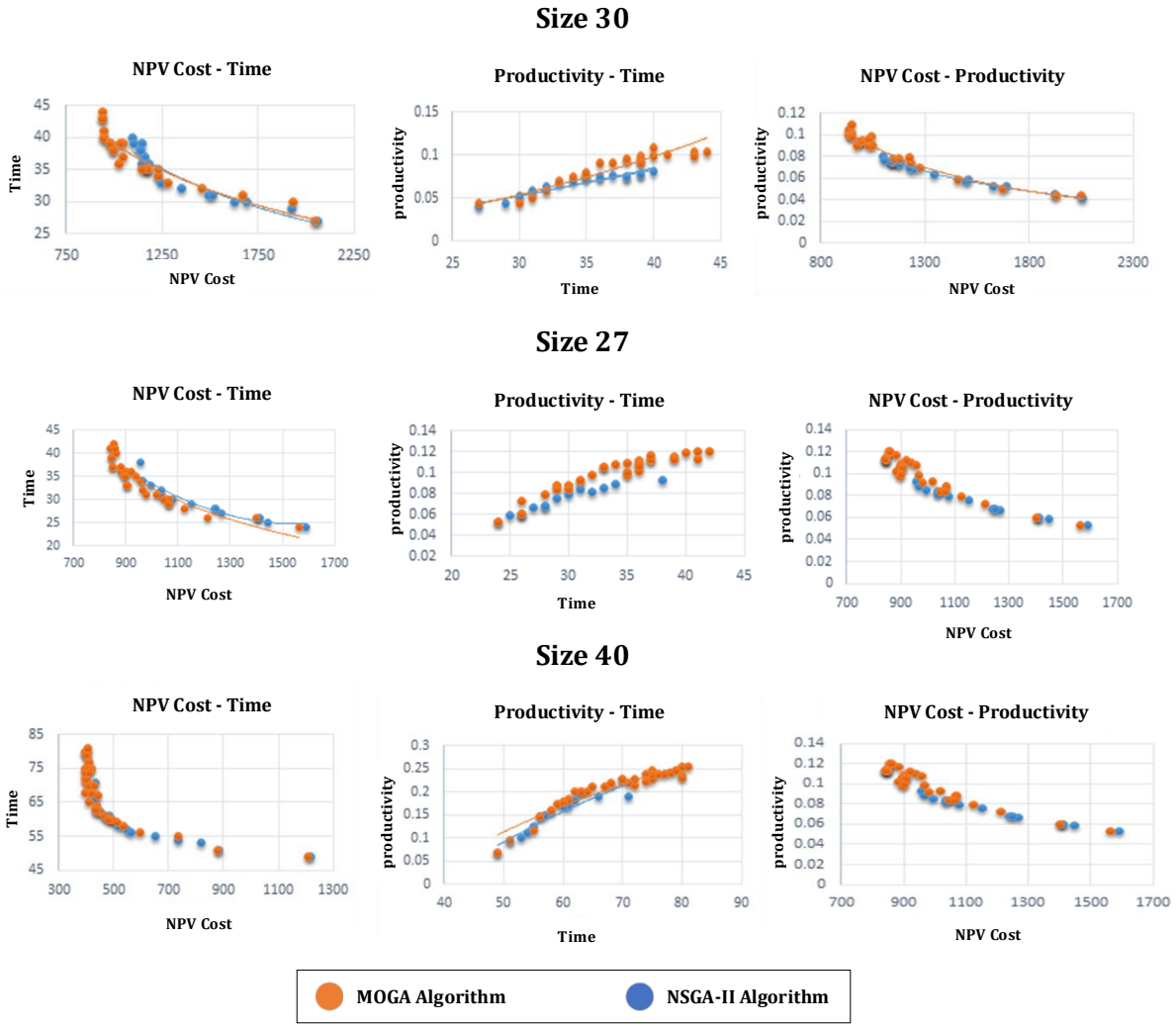

**Figure 16. The Pareto frontier of the selective problems**

In Figure 17 and Figure 18, the solving times of the genetic and NSGA-II algorithms are compared. Evidently, the genetic algorithm takes longer times than NSGA-II to solve all data. In data with fewer activities, the difference between the two algorithms is smaller in their solving times, and this difference strengthens with an increase in the number of points.

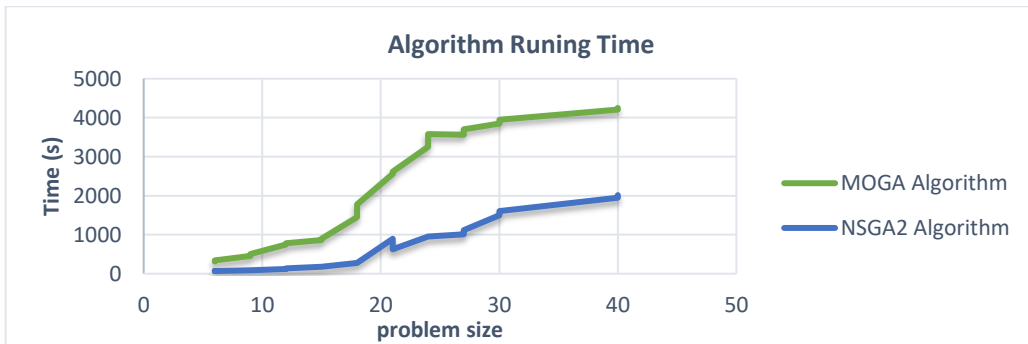

**Figure 17. Diagram comparing the solving time of 23 problems by MOGA and NSGA-II algorithms**

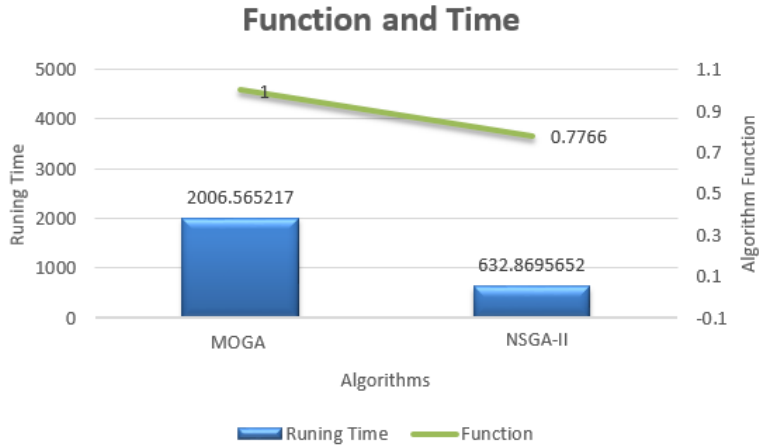

**Figure 18. Diagram comparing the solving time and performance of MOGA and NSGA-II algorithms**

Based on the results of this study and taking into account the NP-hard nature of the problem, it is evident that the accurate method is limited in its ability to find solutions for problems with more than seven activities within a reasonable timeframe. Consequently, the meta-heuristic MOGA algorithm was employed to tackle this challenging problem. The complexity arises primarily from the consideration of multiple executive modes per activity, reflecting the real-world variability in construction projects. These modes enable activities to be completed in shorter timeframes and at lower costs, depending on the chosen mode. When comparing the MOGA algorithm developed in this study with both the accurate method and the classical NSGA-II algorithm, several insights emerge. In small-size problems, MOGA excels by discovering additional Pareto points in significantly less time compared to the accurate solving method. However, in large-size problems, the classic NSGA-II algorithm exhibits superior time performance. Nevertheless, the MOGA algorithm consistently delivers more favorable results according to the specified evaluation criteria, demonstrating its effectiveness. Statistical analysis, specifically the Wilcoxon Signed Ranks Test, confirms the MOGA algorithm's superiority over NSGA-II, with a significant difference in performance. The P-Value calculated for the statistical test is less than 0.05, indicating that the MOGA algorithm outperforms NSGA-II with statistical significance.

## 5.4. Sensitivity Analysis

In this section, the research problem is analyzed based on two essential parameters:

### 5.4.1. Time Horizon Sensitivity Analysis

The impact of time horizon constraints on project completion decisions is substantial and varies depending on the nature of the project. Each project is bound by a specific time horizon within which it must be completed. This time horizon constraint significantly influences the decisions made by contractors. To illustrate this effect, a problem with different time horizons was solved, and the results, as depicted in Figure 19, reveal a clear trend. When the time horizon is shorter, the available choices tend to be less costly but also less productive. Conversely, as the time horizon lengthens, indicating a more extended allowable project completion window, the available choices become more favorable in terms of cost and efficiency. This emphasizes the critical role that time constraints play in project planning and decision-making.

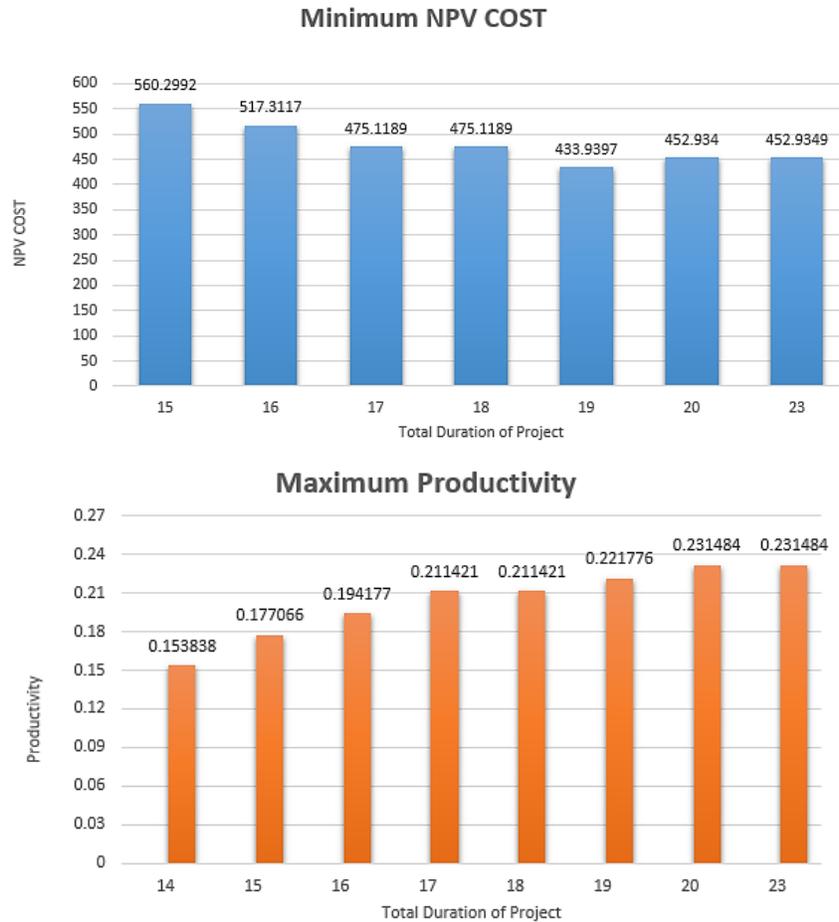

**Figure 19. Time horizon sensitivity analysis**

### 5.4.2. Sensitivity Analysis based on Discount Rate

The discount rate is a fundamental factor in the calculation of the net present value of costs in this research. It represents the interest rate utilized in discounted cash flow (DCF) analysis to ascertain the present value of future cash flows. The discount rate accounts for the time value of money and incorporates considerations of risk or uncertainty associated with future cash flows in cash flow analysis. Specifically, higher uncertainty regarding future cash flows leads to a higher discount rate. To illustrate the impact of this factor on project costs, duration, and productivity, a problem with varying discount rates was analyzed. The results, depicted in Figure 20, reveal a distinct pattern. As the discount rate increases, both the net present value of costs and project completion time shift to the left on the chart. This shift indicates that higher discount rates correspond to lower net present costs (as per the formula for net present costs). Additionally, an increase in the discount rate results in an increase in productivity, as evidenced by its chart. This relationship underscores the intricate interplay between the discount rate, project costs, time, and productivity.

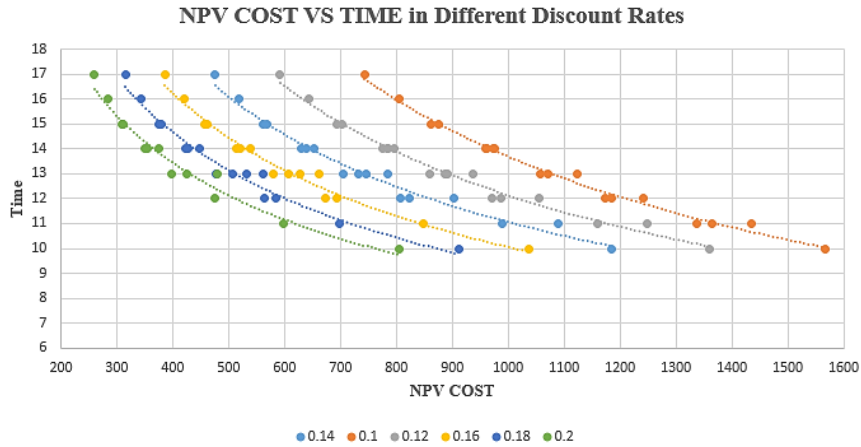

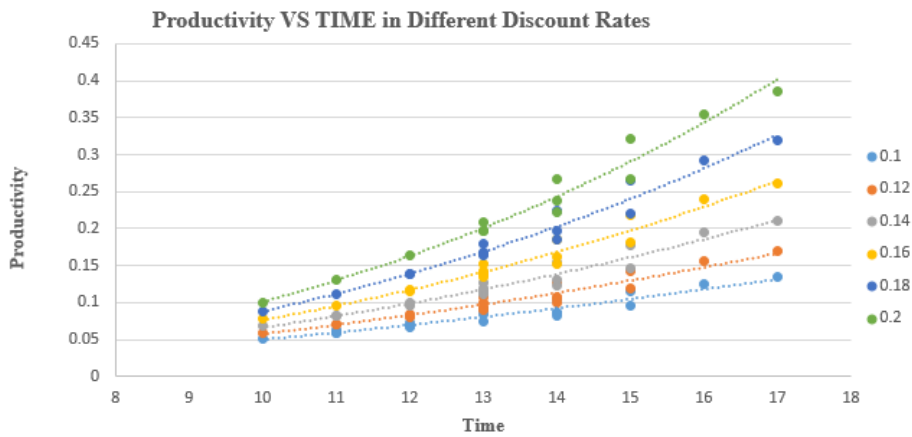

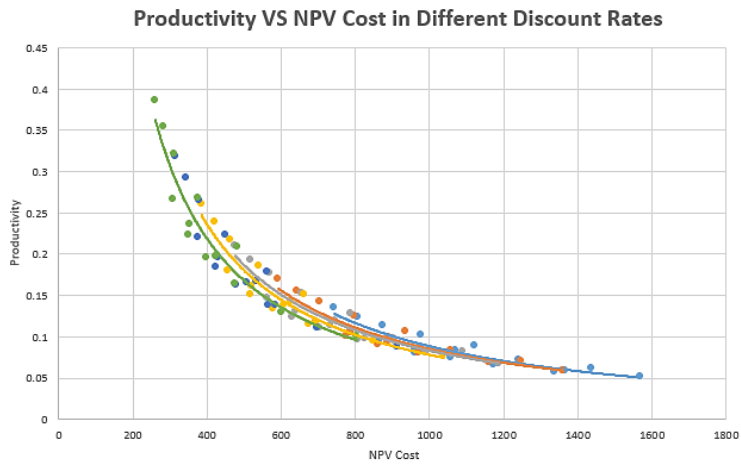

**Figure 20. Sensitivity analysis based on discount rate**

## 6. Conclusion and Future Suggestions

This research delves into the complex realm of multi-mode project scheduling problems, with a specific focus on productivity improvement within construction projects—a topic less explored in the existing literature. The study offers project contractors a comprehensive framework to evaluate and select optimal conditions, balancing minimal delays, enhanced productivity, and cost control. To achieve this, the research formulates the Multimode Resource-Constrained Project Scheduling Problem with Discounted Cash Flow (MMRCPSP-DCF), accounting for varying cost slopes for activities across different execution modes and introducing the Payment at Event Occurrence (PEO) model for payment planning. The resulting three-objective mixed-integer non-linear programming (MINLP) model aims to minimize net present costs, shorten project completion time, and maximize project productivity. Notably, the objectives of cost-time and time-productivity are conflicting, while cost reduction and productivity improvement align. The study employs various solving methods, including the ε-constraint method and two meta-heuristic algorithms (NSGA-II and MOGA), to address different problem sizes and generate optimal Pareto solutions. Performance analysis metrics such as Mean Ideal Distance (MID), Diversity Measure (DM), Number of Pareto Solutions, and Quality Measure (QM) are used to evaluate the solving methods, with MOGA demonstrating effectiveness in both small- and large-size problems. Future research directions could explore uncertainty considerations, environmental impact assessments, and cessation modes for project activities, offering valuable insights for practitioners and decision-makers in project management and construction industries.

# References

Ahmadpour, S. & Ghezavati, V. (2019). Modeling and solving multi-skilled resource-constrained project scheduling problem with calendars in fuzzy condition. *Journal of Industrial Engineering International*, 15, 179-197.

Amoozad Mahdiraji, H., Sedigh, M., Razavi Hajiagha, S.H., Garza-Reyes, J.A., Jafari-Sadeghi, V. & Dana, L.-P. (2021). A novel time, cost, quality and risk tradeoff model with a knowledge-based hesitant fuzzy information: An R&D project application. *Technological Forecasting and Social Change*, 172, 121068.

Audet, C., Bigeon, J., Cartier, D., Le Digabel, S. & Salomon, L. (2021). Performance indicators in multiobjective optimization. *European Journal of Operational Research*, 292, 397-422.

Balouka, N. & Cohen, I. (2021). A robust optimization approach for the multi-mode resource-constrained project scheduling problem. *European Journal of Operational Research*, 291, 457-470.

Birjandi, A. & Mousavi, S.M. (2019). Fuzzy resource-constrained project scheduling with multiple routes: A heuristic solution. *Automation in Construction*, 100, 84-102.

Chen, C.-R., Huang, C.-C. & Tsuei, H.-J. (2014). A Hybrid MCDM Model for Improving GIS-Based Solar Farms Site Selection. *International Journal of Photoenergy*, 2014, 925370.

Chen, P.-H. & Shahandashti, S.M. (2009). Hybrid of genetic algorithm and simulated annealing for multiple project scheduling with multiple resource constraints. *Automation in Construction*, 18, 434-443.

Chen, P.-H. & Weng, H. (2009). A two-phase GA model for resource-constrained project scheduling. *Automation in Construction*, 18, 485-498.

Cui, L., Liu, X., Lu, S. & Jia, Z. (2021). A variable neighborhood search approach for the resource-constrained multi-project collaborative scheduling problem. *Applied Soft Computing*, 107, 107480.

da Silva, A.F., Marins, F.A.S., Tamura, P.M. & Dias, E.X. (2017). Bi-Objective Multiple Criteria Data Envelopment Analysis combined with the Overall Equipment Effectiveness: An application in an automotive company. *Journal of Cleaner Production*, 157, 278-288.

Deb, K., Pratap, A., Agarwal, S., & Meyarivan, T. (2002). A fast and elitist multiobjective genetic algorithm: NSGA-II. *IEEE Transactions On Evolutionary Computation*, 6, 182-197.

Elloumi, S., Loukil, T. & Fortemps, P. (2021). Reactive heuristics for disrupted multi-mode resource-constrained project scheduling problem. *Expert Systems with Applications*, 167, 114132.

Eydi, A. & Bakhshi, M. (2019). A multi-objective resource-constrained project scheduling problem with time lags and fuzzy activity durations. *Journal of Industrial and Systems Engineering*, 12, 45-71.

Fathollahi-Fard, A.M., Hajiaghaei-Keshteli, M. & Tavakkoli-Moghaddam, R. (2020). Red deer algorithm (RDA): a new nature-inspired meta-heuristic. *Soft Computing*, 24, 14637-14665.

Fernandes, G.A. & de Souza, S.R. (2021). A matheuristic approach to the multi-mode resource constrained project scheduling problem. *Computers & Industrial Engineering*, 162, 107592.

Fernandes Muritiba, A.E., Rodrigues, C.D. & Araùjo da Costa, F. (2018). A Path-Relinking algorithm for the multi-mode resource-constrained project scheduling problem. *Computers & Operations Research*, 92, 145-154.

Ghannadpour, S.F., Noori, S., Tavakkoli-Moghaddam, R. & Ghoseiri, K. (2014). A multi-objective dynamic vehicle routing problem with fuzzy time windows: Model, solution and application. *Applied Soft Computing*, 14, 504-527.

Ghannadpour, S.F., & Zandiyeh, F. (2020). A new game-theoretical multi-objective evolutionary approach for cash-in-transit vehicle routing problem with time windows (A Real life Case). *Applied Soft Computing*, 93, 106378.

Ghoddousi, P., Eshtehardian, E., Jooybanpour, S. & Javanmardi, A. (2013). Multi-mode resource-constrained discrete time–cost-resource optimization in project scheduling using non-dominated sorting genetic algorithm. *Automation in Construction*, 30, 216-227.

He, Y., Jia, T. & Zheng, W. (2023). Simulated annealing for centralised resource-constrained multiproject scheduling to minimise the maximal cash flow gap under different payment patterns. *Annals of Operations Research*.

He, Z., He, H., Liu, R. & Wang, N. (2017). Variable neighbourhood search and tabu search for a discrete time/cost trade-off problem to minimize the maximal cash flow gap. *Computers & Operations Research*, 78, 564-577.

Jeunet, J. & Bou Orm, M. (2020). Optimizing temporary work and overtime in the Time Cost Quality Trade-off Problem. *European Journal of Operational Research*, 284, 743-761.

Józefowska, J., Mika, M., Różycki, R., Waligóra, G. & Węglarz, J. (2001). Simulated Annealing for Multi-Mode Resource-Constrained Project Scheduling. *Annals of Operations Research*, 102, 137-155.

Ke, H. & Ma, J. (2014). Modeling project time–cost trade-off in fuzzy random environment. *Applied Soft Computing*, 19, 80-85.

Khalili-Damghani, K., Tavana, M., Abtahi, A.-R. & Santos Arteaga, F.J. (2015). Solving multi-mode time–cost–quality trade-off problems under generalized precedence relations. *Optimization Methods and Software*, 30, 965-1001.

Leyman, P., Van Driessche, N., Vanhoucke, M. & De Causmaecker, P. (2019). The impact of solution representations on heuristic net present value optimization in discrete time/cost trade-off project scheduling with multiple cash flow and payment models. *Computers & Operations Research*, 103, 184-197.

Lin, Y.-K. & Chou, Y.-Y. (2019). A hybrid genetic algorithm for operating room scheduling. *Health Care Management Science*, 1-15.

Liu, H., Fang, Z. & Li, R. (2022). Credibility-based chance-constrained multimode resource-constrained project scheduling problem under fuzzy uncertainty. *Computers & Industrial Engineering*, 171, 108402.

Liu, S.-S. & Wang, C.-J. (2008). Resource-constrained construction project scheduling model for profit maximization considering cash flow. *Automation in Construction*, 17, 966-974.

Liu, Y., Huang, L., Liu, X., Ji, G., Cheng, X. & Onstein, E. (2023). A late-mover genetic algorithm for resource-constrained project-scheduling problems. *Information Sciences*, 642, 119164.

Luong, D.-L., Tran, D.-H. & Nguyen, P.T. (2021). Optimizing multi-mode time-cost-quality trade-off of construction project using opposition multiple objective difference evolution. *International Journal of Construction Management*, 21, 271-283.

Maghsoudlou, H., Afshar-Nadjafi, B. & Akhavan Niaki, S.T. (2017). Multi-skilled project scheduling with level-dependent rework risk; three multi-objective mechanisms based on cuckoo search. *Applied Soft Computing*, 54, 46-61.

Majumdar, J., Verma, M., Shah, P., Karthik, G., Ramachandhran, S. & Gupta, T. (2022). Real-Time Implementation and Analysis of Different Adaptive Enhancement Algorithms Using Embedded Hardware Boards, in: Shetty, N.R., Patnaik, L.M., Nagaraj, H.C., Hamsavath, P.N., Nalini, N. (Eds.), Emerging Research in Computing, Information, Communication and Applications. Singapore: *Springer Singapore*, pp. 1027-1039.

Mehmanchi, E. & Shadrokh, S. (2013). Solving a new mixed integer non-linear programming model of the multi-skilled project scheduling problem considering learning and forgetting effect on the employee efficiency, 2013 IEEE International Conference on Industrial Engineering and Engineering Management, pp. 400-404.

Milička, P., Šůcha, P., Vanhoucke, M. & Maenhout, B. (2022). The bilevel optimisation of a multi-agent project scheduling and staffing problem. *European Journal of Operational Research*, 296, 72-86.

Moosavi Heris, F.S., Ghannadpour, S.F., Bagheri, & M., Zandieh, F. (2022). A new accessibility based team orienteering approach for urban tourism routes optimization (A Real Life Case). *Computers & Operations Research*, 138, 105620.

Mungle, S., Benyoucef, L., Son, Y.-J. & Tiwari, M.K. (2013). A fuzzy clustering-based genetic algorithm approach for time–cost–quality trade-off problems: A case study of highway construction project. *Engineering Applications of Artificial Intelligence*, 26, 1953-1966.

Nabipoor Afruzi, E., Najafi, A.A., Roghanian, E. & Mazinani, M. (2014). A Multi-Objective Imperialist Competitive Algorithm for solving discrete time, cost and quality trade-off problems with mode-identity and resource-constrained situations. *Computers & Operations Research*, 50, 80-96.

Nemati-Lafmejani, R., Davari-Ardakani, H. & Najafzad, H. (2019). Multi-mode resource constrained project scheduling and contractor selection: Mathematical formulation and metaheuristic algorithms. *Applied Soft Computing*, 81, 105533.

Niño, K. & Peña, J. (2019). "A Based-Bee Algorithm Approach for the Multi-Mode Project Scheduling Problem". *Procedia Manufacturing*, 39, 1864-1871.

Panwar, A. & Jha, K.N. (2021). Integrating quality and safety in construction scheduling time-cost trade-off model. *Journal of Construction Engineering and Management*, 147, 04020160.

Pass-Lanneau, A., Bendotti, P. & Brunod-Indrigo, L. (2023). Exact and heuristic methods for Anchor-Robust and Adjustable-Robust RCPSP. *Annals of Operations Research*.

Phuntsho, T. & Gonsalves, T. (2023). Maximizing the Net Present Value of Resource-Constrained Project Scheduling Problems using Recurrent Neural Network with Genetic Algorithm, 2023 International Conference on Intelligent Data Communication Technologies and Internet of Things (IDCIoT): *IEEE*, pp. 524-530.

Rahman, H.F., Chakrabortty, R.K. & Ryan, M.J. (2020). Memetic algorithm for solving resource constrained project scheduling problems. *Automation in Construction*, 111, 103052.

Roghanian, E., Alipour, M. & Rezaei, M. (2018). An improved fuzzy critical chain approach in order to face uncertainty in project scheduling. *International Journal of Construction Management*, 18, 1-13.

Rostami, S., Creemers, S. & Leus, R. (2018). New strategies for stochastic resource-constrained project scheduling. *Journal of Scheduling*, 21, 349-365.

Sadeghi, M.E., Khodabakhsh, M., Ganjipoor, M.R., Kazemipoor, H. & Nozari, H. (2021). A new multi objective mathematical model for relief distribution location at natural disaster response phase. *arXiv preprint arXiv:2108.05458*.

Sadeghloo, M., Emami, S. & Divsalar, A. (2023). A Benders decomposition algorithm for the multi-mode resource-constrained multi-project scheduling problem with uncertainty. *Annals of Operations Research*.

Sayyadi, A., Esmaeeli, H. & Hosseinian, A.H. (2022). A community detection approach for the resource leveling problem in a multi-project scheduling environment. *Computers & Industrial Engineering*, 169, 108202.

Shang, X., Yang, K., Jia, B., Gao, Z. & Ji, H. (2021). Heuristic algorithms for the bi-objective hierarchical multimodal hub location problem in cargo delivery systems. *Applied Mathematical Modelling*, 91, 412-437.

Stylianou, C. & Andreou, A.S. (2016). Investigating the impact of developer productivity, task interdependence type and communication overhead in a multi-objective optimization approach for software project planning. *Advances in Engineering Software*, 98, 79-96.

Su, Y., Lucko, G. & Thompson, R.C. (2020). Apportioning Contract Float with Voting Methods to Correlated Activities in Network Schedules to Protect Construction Projects from Delays. *Automation in Construction*, 118, 103263.

Subulan, K. (2020). An interval-stochastic programming based approach for a fully uncertain multi-objective and multi-mode resource investment project scheduling problem with an application to ERP project implementation. *Expert Systems with Applications*, 149, 113189.

Tareghian, H.R. & Taheri, S.H. (2007). A solution procedure for the discrete time, cost and quality tradeoff problem using electromagnetic scatter search. *Applied Mathematics and Computation*, 190, 1136-1145.

Tavana, M., Abtahi, A.-R. & Khalili-Damghani, K. (2014). A new multi-objective multi-mode model for solving preemptive time–cost–quality trade-off project scheduling problems. *Expert Systems with Applications*, 41, 1830-1846.

Tirkolaee, E.B., Goli, A., Hematian, M., Sangaiah, A.K. & Han, T. (2019). Multi-objective multi-mode resource constrained project scheduling problem using Pareto-based algorithms. *Computing*, 101, 547-570.

Tran, D.-H., Cheng, M.-Y. & Cao, M.-T. (2015). Hybrid multiple objective artificial bee colony with differential evolution for the time–cost–quality tradeoff problem. *Knowledge-Based Systems*, 74, 176-186.

Ulusoy, G., Sivrikaya-Şerifoğlu, F. & Şahin, Ş. (2001). Four Payment Models for the Multi-Mode Resource Constrained Project Scheduling Problem with Discounted Cash Flows. *Annals of Operations Research*, 102, 237-261.

Van Eynde, R., Vanhoucke, M. & Coelho, J. (2023). On the summary measures for the resource-constrained project scheduling problem. *Annals of Operations Research*.

Wang, L. & Zheng, X.-l. (2018). A knowledge-guided multi-objective fruit fly optimization algorithm for the multi-skill resource constrained project scheduling problem. *Swarm and Evolutionary Computation*, 38, 54-63.

Wang, Y., Wang, J.-Q. & Yin, Y. (2020). Multitasking scheduling and due date assignment with deterioration effect and efficiency promotion. *Computers & Industrial Engineering*, 146, 106569.

Yuan, Y., Ye, S., Lin, L. & Gen, M. (2021). Multi-objective multi-mode resource-constrained project scheduling with fuzzy activity durations in prefabricated building construction. *Computers & Industrial Engineering*, 158, 107316.

Zabihi, S., Rashidi Kahag, M., Maghsoudlou, H. & Afshar-Nadjafi, B. (2019). Multi-objective teaching-learning-based meta-heuristic algorithms to solve multi-skilled project scheduling problem. *Computers & Industrial Engineering*, 136, 195-211.